\documentclass[12pt,letterpaper]{article}

\pdfoutput=1

\usepackage{graphics,epsfig,graphicx,color}
\graphicspath{{/EPSF/}{../Figures/}{Figures/}}

\usepackage{amssymb,latexsym}

\usepackage{subfig}
\newsubfloat{figure}

\usepackage{setspace}

\usepackage{amsmath}

\usepackage{mathrsfs,amsfonts}

\usepackage{amsthm,amsxtra}

\newif\ifPDF
\ifx\pdfoutput\undefined
	\PDFfalse
\else
	\ifnum\pdfoutput > 0
		\PDFtrue
	\else
		\PDFfalse
	\fi
\fi

\ifPDF
	\usepackage{pdftricks}
	\begin{psinputs}
		\usepackage{pstricks}
		\usepackage{pstcol}
		\usepackage{pst-plot}
		\usepackage{pst-tree}
		\usepackage{pst-eps}
		\usepackage{multido}
		\usepackage{pst-node}
		\usepackage{pst-eps}
	\end{psinputs}
\else
	\usepackage{pstricks}
\fi

\ifPDF
	\usepackage[debug,pdftex,colorlinks=true, 
	linkcolor=blue, bookmarksopen=false,
	plainpages=false,pdfpagelabels]{hyperref}
\else
	\usepackage[dvips]{hyperref}
\fi

\usepackage{fancyhdr}

\newcommand{\eps}{\epsilon}

\newcommand{\btheta}{{\boldsymbol\theta}}

\newcommand{\bnu}{{\boldsymbol \nu}}

\newcommand{\sfS}{\textsf{S}} \newcommand{\sfE}{\textsf{E}}
\newcommand{\sfC}{\textsf{C}} \newcommand{\sfA}{\textsf{A}}

\newcommand{\bbR}{\mathbb R}

 \newcommand{\bx}{\mathbf x}

 \newcommand{\bJ}{\mathbf J}

 \newcommand{\cT}{\mathcal T}

\newenvironment{keywords}
{\noindent{\bf Key words.}\small}{\par\vspace{1ex}}
\newenvironment{AMS}
{\noindent{\bf AMS subject classifications 2000.}\small}{\par}

\title{On the modeling and simulation of reaction-transfer dynamics in semiconductor-electrolyte solar cells}

\author{
	Yuan He\thanks{Department of Mathematics and ICES, University of Texas, Austin, TX 78712; 
		Email: yuan@ices.utexas.edu .}
	\and 
	Irene M. Gamba\thanks{Department of Mathematics and ICES, University of Texas, Austin, TX 78712; 
		Email: gamba@math.utexas.edu .}
	\and
	Heung-Chan Lee\thanks{Department of Chemistry, University of Texas, Austin, TX 78712; 
		Email: hclee@austin.utexas.edu .}
	\and
   Kui Ren\thanks{Department of Mathematics and ICES, University of Texas, Austin, TX 78712; 
		Email: ren@math.utexas.edu .}
}

\begin{document}

\maketitle



\begin{abstract}
The mathematical modeling and numerical simulation of semiconductor-electrolyte systems play important roles in the design of high-performance semiconductor-liquid junction solar cells. In this work, we propose a macroscopic mathematical model, a system of nonlinear partial differential equations, for the complete description of charge transfer dynamics in such systems. The model consists of a reaction-drift-diffusion-Poisson system that models the transport of electrons and holes in the semiconductor region and an equivalent system that describes the transport of reductants and oxidants, as well as other charged species, in the electrolyte region. The coupling between the semiconductor and the electrolyte is modeled through a set of interfacial reaction and current balance conditions. We present some numerical simulations to illustrate the quantitative behavior of the semiconductor-electrolyte system in both dark and illuminated environments. We show numerically that one can replace the electrolyte region in the system with a Schottky contact only when the bulk reductant-oxidant pair density is extremely high. Otherwise, such replacement gives significantly inaccurate description of the real dynamics of the semiconductor-electrolyte system.
\end{abstract}






\section{Introduction}
\label{SEC:intro}

The mathematical modeling and simulation of semiconductor devices have been extensively studied in past decades due to their importance in industrial applications; see ~\cite{AnAlRi-Book03,Galler-Book05,Grasser-Book03,Jerome-Book96,Jungel-Book09,LaBa-JES76A,LaBa-JES76B,MaRiSc-Book90,Schroeder-Book94,Selberherr-Book84} for overviews of the field and~\cite{Brennan-Book99,HaJa-Book96,Schenk-Book98,YuCa-Book03} for more details on the physics, classical and quantum, of semiconductor devices. In the recent years, the field has been boosted significantly by the increasing need for simulation tools for designing efficient solar cells to harvest sunlight for clean energy. Various theoretical and computational results on traditional semiconductor device modeling are revisited and modified to account for new physics in solar cell applications. We refer interested readers to~\cite{FoPrFeMa-EES12} for a summary of various types of solar cells that have been constructed, to~\cite{KaAtLe-JAP05,Memming-Book01} for simplified analytical solvable models that have been developed, and to~\cite{DePoSaVe-CMAME12,RiPlFoKi-SIAM12,LiChLi-APMM11,Glitzky-M2AS11} for more advanced mathematical and computational analysis of various models. Mathematical modeling and simulation provide ways not only to improve our understanding of the behavior of the solar cells under experimental conditions, but also to predict the performance of solar cells with general device parameters, and thus they enable us to optimize the performance of the cells by selecting the optimal combination of these parameters.

\begin{figure}[!ht]
\begin{minipage}{0.49\textwidth}
\begin{center}
   \includegraphics[angle=0,width=0.98\textwidth]{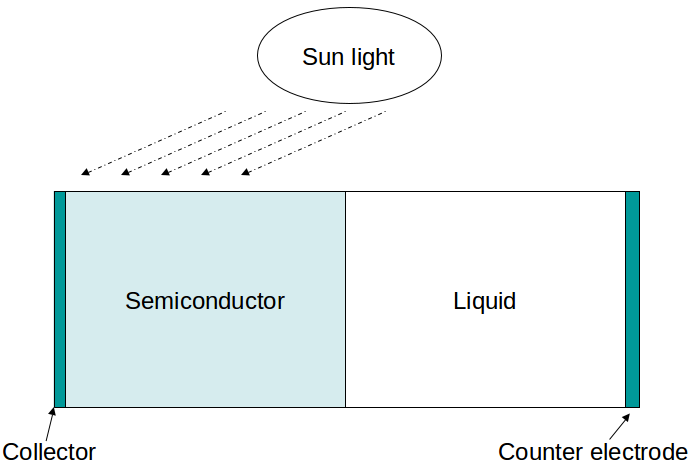}
\end{center}
\end{minipage}
\begin{minipage}{0.25\textwidth}
\vskip 1.4cm
\begin{center}
   \includegraphics[angle=0,width=0.8\textwidth]{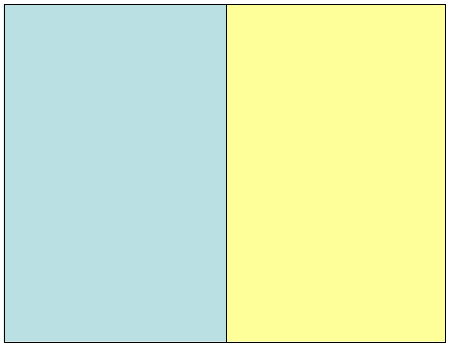}
\end{center}
\end{minipage}
\begin{minipage}{0.24\textwidth}
\vskip 1.5cm
\begin{center}
   \includegraphics[angle=0,width=0.7\textwidth]{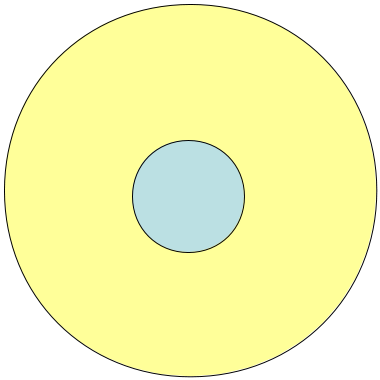}
\end{center}
\end{minipage}
\put(-204,21){$\Gamma_\sfS$}
\put(-204,-62){$\Gamma_\sfS$}
\put(-154,21){$\Gamma_\sfE$}
\put(-154,-62){$\Gamma_\sfE$}
\put(-174,-23){$\Sigma$}
\put(-234,-23){$\Gamma_\sfC$}
\put(-128,-23){$\Gamma_\sfA$}
\put(-45,-25){$\Sigma$}
\put(-15,-25){$\Gamma_\sfA$}
\caption{Left: Sketch of main components in a typical semiconductor-liquid junction solar cell. Middle and Right: Two typical settings for semiconductor-electrolyte systems in dimension two. The semiconductor $\sfS$ and the electrolyte $\sfE$ are separated by the interface $\Sigma$.}
   \label{FIG:Semi-Liq Cell}
\end{figure}
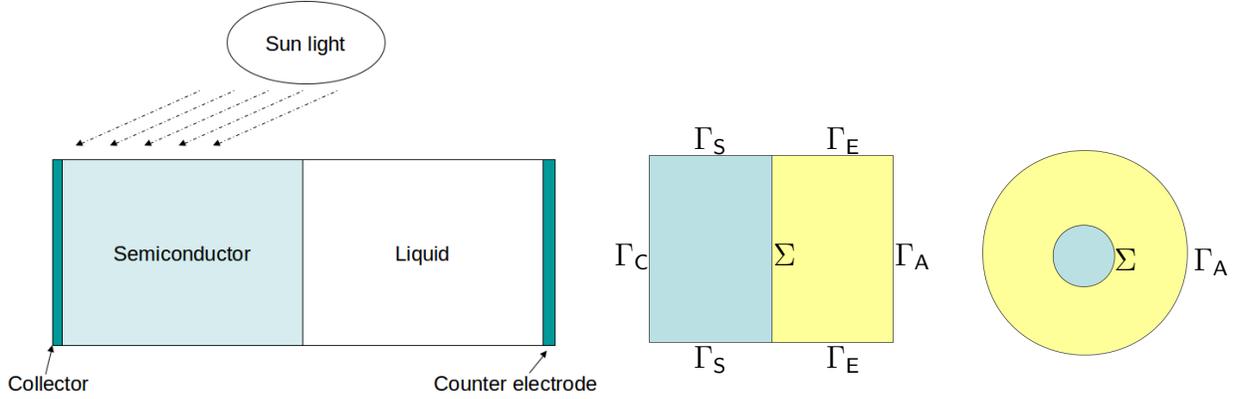

One popular type of solar cells, besides those made of semiconductor p-n junctions, are cells made of semiconductor-liquid junctions. A typical liquid-junction photovoltaic solar cell consists of four major components: the semiconductor, the liquid, the semiconductor-liquid interface and the counter electrode; see a rough sketch in Fig.~\ref{FIG:Semi-Liq Cell} (left). There are many possible semiconductor-liquid combinations; see for instance, ~\cite{FaLe-JPCB97} for ${\rm Si/viologen}^{2+/+}$ junctions, ~\cite{PoLe-JPCB97} for n-type ${\rm InP/Me_2Fc}^{+/0}$ junctions, and ~\cite[Tab. 1]{FoPrFeMa-EES12} for a summary of many other possibilities. The working mechanism of this type of cell is as follows. When sunlight is absorbed by the semiconductor, free conduction electron-hole pairs are generated. These electrons and holes are then separated by an applied potential gradient across the device. The separation of the electrons and holes leads to electrical current in the cell and concentration of charges on the semiconductor-liquid interface where electrochemical reactions and charge transfer occur. We refer interested reader to ~\cite{Gratzel-Nature01,Gratzel-JPPC03,KaTvBaRa-CR10} for physical principles and technical specifics of various types of liquid-junction solar cells.

Charge transport processes in semiconductor-liquid junctions have been studied in the past by many investigators; see~\cite{Lewis-JPCB98} for a recent review. The mechanisms of charge generation, recombination, and transport in both the semiconductor and the liquid are now well understood. However, the reaction and charge transfer process on the semiconductor-liquid interface is far less understood despite the extensive recent investigations from both physical~\cite{GaGeMa-JCP00,GaMa-JCP00,Lewis-ACR90,NoMe-JPC96} and computational~\cite{Memming-Book01,Singh-JAP09} perspectives. The objective of this work is to mathematically model this interfacial charge transfer process so that we could derive a complete system of equations to describe the whole charge transport process in the semiconductor-liquid junction.

To be specific, we consider here semiconductor-liquid junction with the liquid being electrolyte that contains reductant $r$, oxidant $o$ and some other charged species that do not interact with the semiconductor. We denote by $\Omega\subset\bbR^d$ ($d\ge 1$) the domain of interest which contains the semiconductor part $\Omega_{\sfS}$ and the electrolyte part $\Omega_{\sfE}$. We denote by $\Sigma \equiv \partial\Omega_\sfE\cap\partial\Omega_\sfS$ the interface between the semiconductor and electrolyte, $\Gamma_\sfC$ the surface of the current collector at the semiconductor end, and $\Gamma_\sfA$ the surface of the counter (i.e., auxiliary) electrode. $\Gamma_\sfS=\partial\Omega_\sfS\backslash(\Sigma\cup \Gamma_\sfC)$ is the part of the semiconductor boundary that is neither the interface $\Sigma$ nor the contact $\Gamma_\sfC$, and $\Gamma_\sfE=\partial\Omega_\sfE\backslash(\Sigma\cup \Gamma_\sfA)$ the part of electrolyte that is neither the interface $\Sigma$ nor the surface of the counter electrode $\Gamma_\sfA$. We denote by $\bnu(\bx)$ the unit outer normal vector at a point $\bx$ on the boundary of the domain $\partial\Omega=\Gamma_\sfC\cup\Gamma_\sfS\cup\Gamma_\sfA\cup\Gamma_\sfE$. To deal with discontinuities of quantities across the interface $\Sigma$, we use $\Sigma_-$ and $\Sigma_+$ to denote the semiconductor and the electrolyte sides, respectively, of $\Sigma$. On the interface, we use $\bnu^-(\bx)$ and $\bnu^+(\bx)$ to denote the unit normal vectors at $\bx\in\Sigma$ pointing to the semiconductor and the electrolyte domains, respectively.

The rest of the paper is structured as follows. In the next three sections, we present the three main components of the mathematical model. We first model in Section~\ref{SEC:Semi} the dynamics of the electrons and holes in the semiconductor $\Omega_\sfS$. We then model in Section~\ref{SEC:Elec} the dynamics of the reductants and oxidants, as well as other charged species, in the electrolyte $\Omega_\sfE$. In Section~\ref{SEC:Interface} we model the \emph{reaction-transfer dynamics on the interface $\Sigma$}. Once the mathematical model has been constructed, we develop in Section~\ref{SEC:Disc} some numerical schemes for the numerical simulation of the device in simplified settings. We present some numerical experiment in Section~\ref{SEC:Num} where we exhibit the benefits of modeling the complete semiconductor-electrolyte system. Concluding remarks are offered in Section~\ref{SEC:Concl}.

\section{Transport of electrons and holes}
\label{SEC:Semi}

The modeling of transport of free conduction electrons and holes in semiconductor devices has been well studied in recent decades~\cite{AnAlRi-Book03,Galler-Book05,Grasser-Book03,Jerome-Book96,Jungel-Book09,MaRiSc-Book90,Schroeder-Book94,Selberherr-Book84}. Many different models have been proposed, such as the Boltzmann-Poisson system~\cite{BeCaCaVe-JCP09,CaGaMaSh-JCP03,ChGaRe-JCP11,DeNi-NM89,FiJi-JCP10,JiPa-JCP00,Jungel-Book09,MaPi-JCP01,NiDePo-JCP88,Ringhofer-AN97}, the energy transport system~\cite{DrPi-M3AS08,GaJu-SIAM09,Jungel-Book09} and the drift-diffusion-Poisson system~\cite{AnRuRo-SIAM00,CeGaLe-AML97,deMello-JCP02,DeJo-CPC12,Jungel-Book09,MaRiSc-Book90,RiPlFoKi-SIAM12,Rossani-PhysicaA11,SeGaKaAs-SSE07,Steinruck-SIAM89,WaReOd-SIAM91}. For the purpose of computational efficiency, we employ the reaction-drift-diffusion-Poisson model in this work. Let us denote by $(0, T]$ the time interval in which we interested. The bipolar drift-diffusion-Poisson model can be written in the following form:
\begin{equation}\label{EQ:DDP}
  \begin{array}{rcll}
  	\partial_t \rho_n + \nabla\cdot\bJ_n
  	&=& -R_{np}(\rho_n,\rho_p)+\gamma G_{np}(\bx), &\mbox{ in }\ (0, T]\times\Omega_\sfS\\ 
  	\partial_t\rho_p  + \nabla\cdot\bJ_p
  	&=& -R_{np}(\rho_n,\rho_p)+\gamma G_{np}(\bx), &\mbox{ in }\ (0, T]\times\Omega_\sfS\\ 
  	-\nabla\cdot(\eps_r^\sfS\nabla \Phi)&=& 
  	\dfrac{q}{\eps_0}[C(\bx)+\alpha_p\rho_p+\alpha_n\rho_n],
  	&\mbox{ in }\ (0, T]\times\Omega_\sfS. 
  \end{array}
\end{equation}
with the fluxes of electrons and holes are given respectively by 
\begin{equation}\label{EQ:DDP Current}
	\begin{array}{ll}
		\bJ_n = -D_n \nabla \rho_n -\alpha_n \mu_n  \rho_n \nabla \Phi, \qquad 
		\bJ_p = -D_p \nabla \rho_p -\alpha_p \mu_p  \rho_p \nabla \Phi.
	\end{array}
\end{equation}
Here $\rho_n(t,\bx)$ and $\rho_p(t,\bx)$ are the densities of the electrons and the holes, respectively, at time $t$ and location $\bx$, and $\Phi(t,\bx)$ is the electrical potential. The notation $\partial_t$ denotes the derivative with respect to $t$, while $\nabla$ denotes the usual spatial gradient operator. The constant $\eps_0$ is the dielectric constant in vacuum, and the function $\eps_r^\sfS(\bx)$ is the relative dielectric function of the semiconductor material. The function $C(\bx)$ is the doping profile of the device. The coefficients $D_n$ (resp., $D_p$) and $\mu_n$ (resp., $\mu_p$) are the diffusivity and the mobility of electrons (resp., holes). These parameters can be computed from the first principles of statistical physics. In some practical applications, however, they can be fitted from experimental data as well; see, for example, the discussion in~\cite{Schenk-Book98}. The parameter $q$ is the unit electric charge constant, while $\alpha_n=-1$ and $\alpha_p=1$ are the charge numbers of electrons and holes, respectively. The diffusivity and the mobility coefficients are related through the Einstein 
relations $D_n=U_\cT\mu_n$ and $D_p=U_\cT\mu_p$ with $U_\cT$ the thermal voltage at temperature $\cT$ given by $U_\cT=k_B\ \cT/q$, and $k_B$ being the Boltzmann constant.

\subsection{Charge recombination and generation}

The function $R_{np}(\rho_n,\rho_p)$ describes the recombination and generation of electron-hole pairs due to thermal excitation. It represents the rate at which the electron-hole pairs are eliminated through recombination (when $R_{np}>0$) or the rate at which electron-hole pairs are generated (when $R_{np}<0$). Due to the fact that electrons and holes are always recombined and generated in pairs, we have the same rate function for the two species. To be specific, we consider in this work the Auger model of recombination that is based on interactions between multiple electrons and holes, but we refer interested readers to~\cite{Brennan-Book99,FoKiRi-JAP13,KiBe-Book91,MaRiSc-Book90} for discussions on other popular recombination models such as the Shockley-Read-Hall (SRH) model and the Langevin model. 
The Auger model is relevant in cases where the carrier densities are high (for instance in doped materials). It is expressed as
\begin{equation}\label{EQ:Auger}
	R_{np}(\rho_n,\rho_p) = (A_n \rho_n +A_p \rho_p)(\rho_{isc}^2-\rho_n\rho_p),
\end{equation}
where $A_n$ and $A_p$ are the Auger coefficients for electrons and holes respectively. For given materials, $A_n$ and $A_p$ can be measured by experiments. The parameter $\rho_{isc}$ is the intrinsic carrier density that is often calculated from the following formula (see~\cite{KiBe-Book91}):
\begin{equation}\label{EQ:NI}
	\rho_{isc} = \sqrt{N_c N_v}\Big(\dfrac{\cT}{300}\Big)^{1.5} e^{-\frac{E_g}{2k_B \cT}}
\end{equation}
where the band gap $E_g = E_{g0}-\alpha \cT^2/(\cT+\beta)$ with $E_{g0}$ the band gap at $\cT=0 K$ ($E_{g0}=1.17q$ for silicon for instance), $\alpha = 4.73\ 10^{-4}q$, and $\beta = 636$. The parameters $N_c$ and $N_v$ are effective density of states in the conduction and the valence bands, respectively, at $\cT=300 K$.


When the semiconductor device is illuminated by sunlight, the device absorbs photon energy. The absorbed energy creates excitons (bounded electron-hole pairs). The excitons are then separated into free electrons and holes which can then move independently. This generation of free electron-hole pairs is modeled by the source function $G_{np}(\bx)$ in the transport equation~\eqref{EQ:DDP}. Once again, due to the fact that electrons and holes are always generated in pairs, the generating  functions are the same for electrons and holes. We take a model that assumes that photons travel across the device in straight lines. That is, we assume that photons do not get scattered by the semiconductor material during their travel inside the device. This is a reasonable assumption for small devices that have been utilized widely~\cite{KaTvBaRa-CR10}. Precisely, the generation of charges is given as
\begin{equation}\label{EQ:Illu Source}
	G_{np}(\bx)=\left\{
	\begin{array}{rl}
		\sigma(\bx) G_0(\bx_0)e^{-\int_0^{s}\sigma(\bx_0+s'\btheta_0)ds'},& \mbox{if}\ \bx=\bx_0+s\btheta_0\\
	0,& \mbox{otherwise}
	\end{array}\right.
\end{equation}
where $\bx_0\in\Gamma_\sfS$ is the incident location, $\btheta_0$ is the incident direction, $\sigma(\bx)$ is the absorption coefficient (integrated over the usable spectrum), and $G_0(\bx_0)$ is the surface photon flux at $\bx_0$.  The control parameter $\gamma \in\{0,\ 1\}$ in ~\eqref{EQ:DDP} is used to turn on and off the illumination mechanism, and $\gamma=0$ and $\gamma=1$ are the dark and illuminated cases, respectively.

\subsection{Boundary conditions}

We have to supply boundary conditions for the equations in the semiconductor domain. The semiconductor boundary, besides the interface $\Sigma$, is split into two parts, the current collector $\Gamma_\sfC$ and the rest. The boundary condition on the current collector is determined by the type of contacts formed there. There are mainly two types of contacts, the Ohmic contact and the Schottky contact.

\paragraph{Dirichlet at Ohmic contacts.} Ohmic contacts are generally used to model metal-semiconductor junctions that do not rectify current. They are appropriate when the Fermi levels in the metal contact and adjacent semiconductor are aligned. Such contacts are mainly used to carry electrical current out and into semiconductor devices, and should be fabricated with little (or ideally no) parasitic resistance. Low resistivity Ohmic contacts are also essential for high-frequency operation. Mathematically, Ohmic contacts are modeled by Dirichlet boundary conditions which can be written as~\cite{Brennan-Book99,KiBe-Book91,MaRiSc-Book90,Nelson-Book03}
\begin{equation}\label{EQ:DDP BC Dirich}
	\begin{array}{lccl}
	\rho_n(t,\bx)=\rho_n^e(\bx), & & \rho_p(t,\bx)=\rho_p^e(\bx), & \mbox{on}\ \ (0, T]\times\Gamma_{\sfC},\\
  	\Phi(t,\bx) =\varphi_{bi}+\varphi_{app}, & & & \mbox{on}\ \ (0, T]\times\Gamma_{\sfC},
  \end{array}
\end{equation}
where $\varphi_{bi}$ and $\varphi_{app}$ are the built-in and applied potentials, respectively. The boundary values $\rho_n^e$, $\rho_p^e$ for the Ohmic contacts are calculated following the assumptions that the semiconductor is in stationary and equilibrium state and that the charge neutrality condition holds. This means that right-hand-side of the Poisson equation disappears so that 
\begin{equation}\label{EQ:Equi}
	C+\rho_p^e-\rho_n^e=0.
\end{equation}
Thermal equilibrium implies that generation and recombination balance out, so $R_{np}=0$ at Ohmic contacts. This leads to the mass-action law, between the density of electrons and holes: 
\begin{equation}\label{EQ:Mass}
	\rho_n^e\rho_p^e-\rho_{isc}^2=0.
\end{equation}
The system of equations~\eqref{EQ:Equi} and ~\eqref{EQ:Mass} admit a unique solution pair $(\rho_n,\rho_p)$, which is given by
\begin{equation}\label{EQ:Ohmic}
	\begin{array}{l}
		\rho_n^e(t,\bx)=\dfrac{1}{2}(\sqrt{C^2+4\rho_{isc}^2}+C),\qquad 		
		\rho_p^e(t,\bx)=\dfrac{1}{2}(\sqrt{C^2+4\rho_{isc}^2}-C) .
	\end{array}
\end{equation}
These densities result in a built-in potential that can be calculated as
\begin{equation}
	\varphi_{bi} = U_\cT \ln(\rho_n^e/\rho_{isc}).
\end{equation}
Note that due to the fact that the doping profile $C$ varies in space, these boundary values are different on different part of the boundary.

\paragraph{Robin (or mixed) at Schottky contacts.} Schottky contacts are used to model metal-semiconductor junctions that have rectifying effects (in the sense that current flow through the contacts is rectified). They are appropriate for contacts between a metal and a lightly doped semiconductor. Mathematically, at a Schottky contact, Robin- (also called mixed-) type of boundary conditions are imposed for the ${\rm n}$- and ${\rm p}$-components, while Dirichlet-type of conditions are imposed for the $\Phi$-component. More precisely, these boundary conditions are~\cite{Brennan-Book99,KiBe-Book91,MaRiSc-Book90,Nelson-Book03}:
\begin{equation}\label{EQ:DDP BC Robin}
	\begin{array}{lccl}
	\bnu\cdot\bJ_n(t,\bx)=v_n(\rho_n-\rho_n^e)(\bx), & & 
	\bnu\cdot\bJ_p(t,\bx)=v_p(\rho_p-\rho_p^e)(\bx), & \mbox{on}\ \ (0, T]\times\Gamma_{\sfC},\\
  	\Phi(t,\bx) =\varphi_{Stky}+\varphi_{app}, & &  & \mbox{on}\ \ (0, T]\times\Gamma_{\sfC}.
  \end{array}
\end{equation}
Here the parameters for the Schottky barrier are the recombination velocities $v_n$ and $v_p$, and the height of the potential barrier, $\varphi_{Stky}$, which depends on the materials of the semiconductor and the metal in the following way:
\begin{equation}\label{EQ:Schottky Height}
	\varphi_{Stky} = \left\{
	\begin{array}{ll}
		\Phi_m-\chi, & \mbox{n-type}\\
		\frac{E_g}{q}-(\Phi_m-\chi), & \mbox{p-type}
	\end{array}
	\right.
\end{equation}
where $\Phi_m$ is the work function, i.e., the potential difference between the Fermi energy and the vacuum level, of the metal and $\chi$ is the electron affinity, i.e., the potential difference between the conduction band edge and the vacuum level. $E_g$ is again the band gap. The values of the parameters $v_n,\ v_p,\ \Phi_m$, and $\chi$ are given in Tab.~\ref{TAB:Para} of Section~\ref{SEC:Num}.

\paragraph{Neumann at insulating boundaries.} On the part of the semiconductor boundary that is not the current collector, it is natural to impose insulating boundary conditions which ensures that there is no charge or electrical currents through the boundary. The conditions are
\begin{equation}\label{EQ:DDP BC Neumann}
	\begin{array}{lccl}
	\bnu\cdot D_n\nabla\rho_n(t,\bx)=0, & & 
	\bnu\cdot D_p\nabla\rho_p(t,\bx)=0, & \mbox{on}\ \ (0, T]\times\Gamma_{\sfS},\\
  	\bnu\cdot \eps_r^\sfS \nabla\Phi(t,\bx) =0, & &  & \mbox{on}\ \ (0, T]\times\Gamma_{\sfS}.
  \end{array}
\end{equation}
In solar cell applications, part of the boundary $\Gamma_\sfS$ is where illumination light enters the semiconductor.

We finish this section with the following remarks. It is generally believed that the Boltzmann-Poisson model~\cite{Jungel-Book09} is a more accurate model for charges transport in semiconductors. However, the Boltzmann-Poisson model is computationally more expensive to solve and analytically more complicated to analyze. The drift-diffusion-Poisson model~\eqref{EQ:DDP} can be regarded as a macroscopic approximation to the Boltzmann-Poisson model. The validity of the drift-diffusion-Poisson model can be justified in the case when the mean free path of the charges is very small compared to the size of the device and the potential drop across the device is small (so that the electric field is not strong); see, for instance, ~\cite{BeDe-JMP96,CaGaSh-PhysicaD00,CeGaLe-AML97,CeGaLe-SIAM01,Jungel-Book09} for such a justification.

\section{Charge transport in electrolytes}
\label{SEC:Elec}

We now present the equations for the reaction-transport dynamics of charges in an electrolyte. To be specific, we consider here only electrolytes that contain reductant-oxidant pairs (denoted by $r$ and $o$) that interact directly with the semiconductor through electrons transfer (which we will model in the next section), and $N$ other charged species (denoted by $j=1,...,N$) that do not interact directly with the semiconductor through electron transfer. We also limit our modeling efforts to reaction, recombination, transport, and diffusion of the charges. Other more complicated physical and chemical processes are neglected.

We model the charge transport dynamics in electrolyte again with a set of reaction-drift-diffusion-Poisson equations. In the electrochemistry literature, this mathematical description of the dynamics is often called the Poisson-Nernst-Planck theory~\cite{BaThAj-PRE04,BiFuBa-RJE12,Eisenberg-SM00,Fawcett-Book04,HoLiLiEi-JPCB12,Liu-JDE09,KiBaAj-PRE07B,LuZh-BJ11,MaMu-IJHMT09,MaPeAg-JCP88,PaJe-SIAM97,ScNaEi-PRE01,SiNo-SIAM09}. Let us denote by $\rho_r(t,\bx)$ the density of the reductants, by $\rho_o(t,\bx)$ the density of the oxidants, and by $\rho_j$ ($1\le j\le N$) the density of the other $N$ charge species. Then these densities solve the following system that is of the same form as~\eqref{EQ:DDP}:
\begin{equation}\label{EQ:DDP Redox}
	\begin{array}{ll}
		\partial_t \rho_r + \nabla\cdot \bJ_r = +R_{ro}(\rho_r,\rho_o), 
		&\text{in }\  (0, T]\times\Omega_\sfE,\\
		\partial_t \rho_o + \nabla\cdot \bJ_o = -R_{ro}(\rho_r,\rho_o), 
		&\text{in }\  (0, T]\times\Omega_\sfE,\\
		\partial_t \rho_j + \nabla\cdot \bJ_j = R_{j}(\rho_1,\cdots,\rho_N),\ \ 1\le j\le N
		&\text{in }\  (0, T]\times\Omega_\sfE,\\ 
		-\nabla\cdot\eps_r^\sfE\nabla\Phi= \dfrac{q}{\eps_0} (\alpha_o\rho_o+\alpha_r\rho_r+\sum_{j=1}^N\alpha_j\rho_j), 
		&\text{in }\  (0, T]\times\Omega_\sfE.
\end{array}
\end{equation}
with the fluxes given respectively by
\begin{equation}\label{EQ:DDP Redox Current}
\begin{array}{lcl}
	\bJ_r = -D_r \nabla\rho_r - \alpha_r\mu_r \rho_r \nabla \Phi, & & 
	\bJ_o = -D_o \nabla\rho_o - \alpha_o\mu_o \rho_o \nabla \Phi\\
	\bJ_j = -D_j \nabla\rho_j - \alpha_j\mu_j \rho_j \nabla \Phi, & & 1\le j\le N
\end{array}
\end{equation}
where again the diffusion coefficient $D_r$ (resp., $D_o$ and $D_j$) is related to the mobility $\mu_r$ (resp., $\mu_o$ and $\mu_j$) through the Einstein relation $D_r=U_\cT \mu_r$ (resp., $D_o=U_\cT \mu_o$ and $D_j=U_\cT \mu_j$). The parameters $\alpha_o$, $\alpha_r$ and $\alpha_j$ ($1\le j\le N$) are the charge numbers of the corresponding charges species. Depending on the type of the redox pairs in the electrolyte, the charge numbers can be different; see, for instance, ~\cite{FoPrFeMa-EES12} for a summary of various types of redox electrolytes that have been developed.

Let us remark that in the above modeling of the dynamics of reductant-oxidant pair, we have implicitly assumed that the electrolyte, in which the redox pairs live, is not perturbed by charge motions. In other words, there is no macroscopic deformation of the electrolyte that can occur. If this is not the case, we have to introduce the equations of fluid dynamics, mainly the Navier-Stokes equation, for the fluid motion, and add an advection term (with advection velocity given by the solution of the Navier-Stokes equation) in the current expressions in~\eqref{EQ:DDP Redox Current}. The dynamics will thus be far more complicated.

\subsection{Charge generation through reaction}

The reaction mechanism between the oxidants and the reductants is modeled by the reaction rate function $R_{ro}$. Note that the elimination and generation of the redox pairs are different from those of the electrons and holes. An oxidant is eliminated (resp., generated) when a reductant is generated (resp., eliminated) and vice versa. This is the reason why there is a negative sign in front of the function $R_{ro}$ in the second equation of ~\eqref{EQ:DDP Redox}. The oxidation-reduction reaction requires free electrons which are only available through the semiconductor. Therefore this reaction occurs mainly on the semiconductor-electrolyte interface. We thus assume in general that there is no oxidation-reduction reaction in the bulk electrolyte; that is,
\begin{equation}
	R_{ro}(\rho_r,\rho_o)=0,\quad \mbox{in}\ \ (0, T]\times\Omega_\sfE .
\end{equation}
This is what we adopt in the simulations of Section~\ref{SEC:Num}. The oxidation-reduction reaction on the interface $\Sigma$ will be modeled in the next section.

The reactions among other charged species in the electrolyte are modeled by the reaction rate functions $R_j$ ($1\le j\le N$). The exact forms of these rate functions can be derived following the law of mass action once the types of reactions among the charged species presented in the electrolyte are known. We refer to ~\cite{CaBi-Book72,Connors-Book90} for the rate functions for various chemical reactions.

\subsection{Boundary conditions}

It is generally assumed that the interface of semiconductor and electrolyte is far from the counter electrode. Therefore the values for the densities of redox pairs on the electrode $\Gamma_\sfA$ are set as their bulk values. Mathematically, this means that Dirichlet boundary conditions have to be imposed:
\begin{equation}\label{EQ:DDP Redox BC}
	\begin{array}{lll}
	\rho_r(t,\bx)=\rho_r^\infty(\bx), & \rho_o(t,\bx)=\rho_o^\infty(\bx), & \mbox{on}\ \ (0, T]\times\Gamma_{\sfA},\\
  	\rho_j(t,\bx)=\rho_j^\infty(\bx),\ 1\le j\le N, & \Phi(t,\bx) =\varphi_{app}^\sfA, & \mbox{on}\ \ (0, T]\times\Gamma_{\sfA},
  \end{array}
\end{equation}
where $\rho_r^\infty$ and $\rho_o^\infty$ are the bulk concentration of the respective species, and  $\varphi_{app}^\sfA$ is the applied potential on the counter electrode. The values of these parameters are given in Tab.~\ref{TAB:Para} in Section~\ref{SEC:Num}.

On the rest of the electrolyte boundary, $\Gamma_\sfE$, we impose again insulating boundary conditions:
\begin{equation}\label{EQ:DDP Redox BC Neumann}
	\begin{array}{lll}
	\bnu\cdot D_r\nabla \rho_r(t,\bx)=0, & \bnu\cdot D_o\nabla\rho_o(t,\bx)=0, & \mbox{on}\ \ (0, T]\times\Gamma_{\sfE},\\
  	\bnu\cdot D_j\nabla \rho_j(t,\bx)=0,\ 1\le j\le N, & \bnu\cdot \eps_r^\sfE\nabla\Phi(t,\bx) =0, & \mbox{on}\ \ (0, T]\times\Gamma_{\sfE},
  \end{array}
\end{equation}

\section{Interfacial reaction and charge transfer}
\label{SEC:Interface}

In order to obtain a complete mathematical model for the semiconductor-electrolyte system, we have to couple the system of equations in the semiconductor with those in the electrolyte through interface conditions that describe the interfacial charge transfer process.

\subsection{Electron transfer between electron-hole and redox}

The microscopic electrochemical processes on the semiconductor-electrolyte interface can be very complicated, depending on the types of semiconductor materials and electrolyte solutions. There is a vast literature in physics and chemistry devoted to the subject; see, for instance, ~\cite{BeFaPeWi-EMAC03,GaGeMa-JCP00,GaMa-JCP00,KaTvBaRa-CR10,Lewis-ACR90,Lewis-JPCB98,Memming-Book01,NoMe-JPC96,PeFaPl-SEMSC08,PeFaWiBe-JPP04,Singh-JAP09} and references therein. We are only interested in deriving macroscopic interface conditions that are consistent with the dynamics of charge transport in the semiconductor and the electrolyte modeled by the equation systems~\eqref{EQ:DDP} and ~\eqref{EQ:DDP Redox}. Without attempting to construct models in the most general cases, we focus here on oxidation-reduction reactions described by the following process,
\begin{equation}
	\mbox{Ox}^{|\alpha_0|+} + e^{-}(\mbox{\sfS}) \rightleftharpoons \mbox{Red}^{|\alpha_r|-},
\end{equation}
with $\alpha_o-\alpha_r=1$, $\mbox{Ox}$ and $\mbox{Red}$ denoting respectively the oxidant and the reductant, and $e^-(\sfS)$ denoting an electron from the semiconductor. Experimental studies semiconductor-electrolyte interface with this type of reaction can be found in~\cite{FaLe-JPCB97} for ${\rm Si/viologen}^{2+/+}$ interfaces and in~\cite{PoLe-JPCB97} for n-type ${\rm InP/Me_2Fc}^{+/0}$ interfaces.

The changes of the concentrations of the redox pairs on $\Sigma_+$, after taking into account the conservation of $\rho_r+\rho_o$, can be written respectively as:
\begin{equation}
	\dfrac{d{\rho_r}}{dt}={k_f}{\rho_o}-{k_b}{\rho_r}
	\qquad \mbox{and} \qquad  
	\dfrac{d{\rho_o}}{dt}={k_b}{\rho_r}-{k_f}{\rho_o}, 
\end{equation}
where $k_f$ and $k_b$ are the pseudo first-order forward and backward reaction rates, respectively. The forward reaction rate $k_f$ is proportional to the product of the electron transfer rate $k_{et}$ through the interface and the density of the electrons on $\Sigma_-$. The backward reaction rate $k_b$ is proportional to the product of the hole transfer rate $k_{ht}$ and the density of the holes on $\Sigma_-$. More precisely, we have~\cite{GaGeMa-JCP00,GaMa-JCP00,Lewis-JPCB98,Memming-Book01}:
\begin{equation}
	k_f(t,\bx) =k_{et}(\bx)(\rho_n-\rho_n^e) \ \ \mbox{and} \ \ 
	k_b(t,\bx) =k_{ht}(\bx)(\rho_p-\rho_p^e),\ \ \mbox{on}\ \ (0,T]\times\Sigma .
\end{equation}

The changes of the concentrations lead to, following the relations $\bnu^+\cdot\bJ_r =-\frac{d{\rho_r}}{dt}$ and $\bnu^+\cdot\bJ_o =-\frac{d{\rho_o}}{dt} $, fluxes of redox pairs through the interface that can be expressed as follows~\cite{GaGeMa-JCP00,GaMa-JCP00,Lewis-JPCB98,Memming-Book01}:
\begin{equation}\label{EQ:Interf E}
	\begin{array}{ll}
	\bnu^+\cdot\bJ_r ={k_{ht}(\bx)(\rho_p-\rho_p^e)}{\rho_r(t,\bx)}-{k_{et}(\bx)(\rho_n-\rho_n^e)}{\rho_o}(t,\bx), & \mbox{on}\ \ (0,T]\times\Sigma\\
	\bnu^+\cdot\bJ_o =-{k_{ht}(\bx)(\rho_p-\rho_p^e)}{\rho_r(t,\bx)}+{k_{et}(\bx)(\rho_n-\rho_n^e)}{\rho_o(t,\bx)}, & \mbox{on}\ \ (0,T]\times\Sigma
	\end{array}
\end{equation}
where the unit normal vector $\bnu^+$ points toward the electrolyte domain.

The fluxes of the redox pairs from the interface given in~\eqref{EQ:Interf E} consists of two contributions: the flux induced from the transfer of electrons from the semiconductor to the electrolyte, often called the cathodic current after being brought up to the right dimension, $\bnu^+\cdot\bJ_n$, and the flux induced from the transfer of electrons from the electrolyte to the conduction band, often called the anodic current  after being brought up to the right dimension, $\bnu^+\cdot\bJ_p$:
\begin{equation}\label{EQ:Interf S}
	\begin{array}{rcll}
	\bnu^+\cdot\bJ_n=(-\bnu^-)\cdot\bJ_n &=& - k_{et}(\bx)(\rho_n-\rho_n^e) \rho_o(t,\bx),& \mbox{on}\ \ (0,T]\times\Sigma\\ 
	\bnu^+\cdot\bJ_p=(-\bnu^-)\cdot\bJ_p &=& - k_{ht}(\bx)(\rho_p-\rho_p^e)\rho_r(t,\bx), & \mbox{on}\ \ (0,T]\times\Sigma
	\end{array}
\end{equation}

The interface conditions~\eqref{EQ:Interf E} and ~\eqref{EQ:Interf S} can now be supplied to the semiconductor equations in~\eqref{EQ:DDP} and the redox equations in ~\eqref{EQ:DDP Redox} respectively. The values of the electron and hole transfer rates, $k_{et}$ and $k_{ht}$, in the interface conditions can be calculated approximately from the first principles of physical chemistry~\cite{BeFaPeWi-EMAC03,GaGeMa-JCP00,GaMa-JCP00,KaTvBaRa-CR10,Lewis-ACR90,Memming-Book01,NoMe-JPC96,PeFaPl-SEMSC08,PeFaWiBe-JPP04,Singh-JAP09}. Theoretical analysis shows that both parameters can be approximately treated as constant; see, for instance, ~\cite{Lewis-JPCB98} for a summary of various ways to approximate these rates. The dependence of the forward and reverse reaction rates, $k_f$ and $k_b$, on the electric potential $\Phi$ is encoded in the electron and hole densities (i.e. $\rho_n$ and $\rho_p$) on the interface. In fact, we can recover the commonly used Butler-Volmer model~\cite{BaFa-Book00,NeTh-Book04} from our model as follows. Consider a one-dimensional semiconductor-electrolyte system (just for the purposes of presentation). At the equilibrium of the system, the net reaction rate is zero, i.e $\bnu^+\cdot\bJ_o=-\bnu^+\cdot\bJ_r=\bnu^+\cdot \bJ_n-\bnu^+\cdot \bJ_p=k_f\rho_o-k_b\rho_r=0$ on $\Sigma_+$. This leads to, by the expression for the fluxes~\eqref{EQ:DDP Redox Current}, the following relation between densities of redox pairs and the electric potential: 
\begin{equation}
	\rho_r=\rho_r^e \exp\left(\frac{\Phi-\Phi^e}{U_\cT}\right), \qquad 
	\rho_o=\rho_o^e \exp\left(\frac{\Phi-\Phi^e}{U_\cT}\right)
\end{equation}
where we have used the Einstein relations $D_r=U_\cT\mu_r$ and $D_o=U_\cT\mu_o$, $\Phi^e$ is the equilibrium potential of the electrolyte, and $\rho_r^e$, $\rho_o^e$ are the corresponding densities. Therefore at the system equilibrium the forward and reverse reaction rates satisfy:
\begin{equation}
	\frac{k_f}{k_b}=\dfrac{\rho_o}{\rho_r}=\frac{\rho_o^e}{\rho_r^e} \exp\left(\frac{-2(\Phi-\Phi^e)}{U_\cT}\right).
\end{equation}
This implies that
\begin{equation}
	-\dfrac{1}{2} U_\cT \left(\dfrac{d \ln k_f}{d\Phi}-\dfrac{d \ln k_b}{d\Phi}\right)=1.
\end{equation}
If we define the reductive symmetry factor (related to the forward reaction) $\alpha= -\frac{U_\cT}{2} \frac{d \ln k_f}{d\Phi}$, the above relation then implies that $\frac{U_\cT}{2} \frac{d \ln k_b}{d\Phi}=1-\alpha$. Moreover, the definition leads to $k_f=k^0\exp(-\alpha\eta(\Phi-\Phi^e))$ and $k_b=k^0\exp((1-\alpha)\eta (\Phi-\Phi^e))$ with $\eta=2/U_\cT$, where $k^0$ is the standard rate constant. We can therefore have the following Butler-Volmer model~\cite{BaFa-Book00,JiWaSiKaLeXi-EES14,NeTh-Book04} for the interface flux:
\begin{equation}\label{EQ:BV}
	\bnu^+\cdot\bJ_o = -\bnu^+\cdot\bJ_r =k^0\big[e^{-\alpha\eta(\Phi-\Phi^e)}{\rho_o}(t,\bx) -e^{(1-\alpha)\eta (\Phi-\Phi^e)} \rho_r(t,\bx)\big].
\end{equation}

\subsection{Interface conditions for nonredox charges}

For the $N$ charged species in the electrolyte that do not interact directly with the semiconductor through electron transfer, we impose insulating boundary conditions on the interface:
\begin{equation}\label{EQ:Interf Non Redox}
	\begin{array}{lll}
	\bnu^-\cdot \bJ_j = 0,&   1\le j\le N, & \mbox{on}\ \ (0, T]\times\Sigma
	\end{array}
\end{equation}

We need to specify the interface condition for the electric potential as well. This is done by requiring $\Phi$ to be continuous across the interface and have continuous flux. Let us denote by $\Sigma_+$ and $\Sigma_-$ the semiconductor and the electrolyte sides, respectively, of $\Sigma$; then the conditions on the electric potential are given by
\begin{equation}\label{EQ:Interf Phi}
	[\Phi]_{\Sigma}\equiv\Phi_{|\Sigma_-}-\Phi_{|\Sigma_+}=0, \qquad 
	[\eps_r\frac{\partial\Phi}{\partial\nu}]_{\Sigma}\equiv \Big(\eps_r^\sfE\dfrac{\partial\Phi}{\partial\nu}\Big)_{|\Sigma_-}-\Big(\eps_r^\sfS\dfrac{\partial\Phi}{\partial\nu}\Big)_{|\Sigma_+}=0, \quad \mbox{on}\ \ (0, T]\times\Sigma.
\end{equation}
Note that these continuity conditions would \emph{not} prevent large the electric potential drops across a narrow neighborhood of the interface, as we will see in the numerical simulations.

The interface conditions that we constructed in this section ensure the conservation of the total flux $\bJ$ across the interface. To check that we recall that the total flux in the system is given as~\cite{CaBi-Book72,Fawcett-Book04,Hille-Book01,Selberherr-Book84}
\begin{equation}\label{EQ:Total Current}
	\bJ(\bx)=\left\{
	\begin{array}{cl}
		\alpha_p\bJ_p+\alpha_n\bJ_n, & \bx\in\sfS\\
		\alpha_o\bJ_o+\alpha_r\bJ_r+\sum_{j=1}^N\alpha_j\bJ_j, & \bx\in\sfE
	\end{array}\right.	
\end{equation}
Using~\eqref{EQ:Interf E},~\eqref{EQ:Interf S}, ~\eqref{EQ:Interf Non Redox} and the fact that $\alpha_o-\alpha_r=1$, we verify that $\bnu^+\cdot\bJ_{|\Sigma_-}=\bnu^+\cdot\bJ_{|\Sigma_+}$.

\section{Numerical discretization}
\label{SEC:Disc}

The drift-diffusion-Poisson equations in ~\eqref{EQ:DDP} and ~\eqref{EQ:DDP Redox}, together with the boundary conditions given in~\eqref{EQ:DDP BC Dirich} (resp., ~\eqref{EQ:DDP BC Robin}), ~\eqref{EQ:DDP BC Neumann}, ~\eqref{EQ:DDP Redox BC}, ~\eqref{EQ:DDP Redox BC Neumann}, and the interface conditions~\eqref{EQ:Interf E}, ~\eqref{EQ:Interf S}, ~\eqref{EQ:Interf Non Redox}, and ~\eqref{EQ:Interf Phi} form a complete mathematical model for the transport of charges in the system of semiconductor-electrolyte for solar cell simulations. We now present a numerical procedure to solve the system. 

\subsection{Nondimensionalization}

We first introduce the following characteristic quantities in the simulation regarding the device and its physics. We denote by $l^*$ the characteristic length scale of the device, $t^*$ the characteristic time scale, $\Phi^*$ the characteristic voltage, and $C^*$ the characteristic density. The values (and units) for these characteristic quantities are respectively (see ~\cite{KiBe-Book91,LaBa-JES76A,LaBa-JES76B,Nelson-Book03,PeFaPl-SEMSC08}),
\begin{equation}\label{EQ:Characteristic}
l^*=10^{-4}\ {\rm (cm)},\ \ \  t^*=10^{-12}\ {\rm (s)},\ \ \ \Phi^*=U_\cT\ {\rm (V)},\ \ \ C^*= 10^{16}\ {\rm (cm^{-3})} .
\end{equation}

We now rescale all the physical quantities. For any quantity $Q$, we use $Q'$ to denote its rescaled version. To be specific, we introduce the rescaled Debye lengths in the semiconductor and electrolyte regions, respectively, as
\begin{displaymath}
\lambda_\sfS=\frac{1}{l^*}\sqrt{\frac{\Phi^*\epsilon^\sfS}{q C^*}},\quad \mbox{and}\quad \lambda_\sfE=\frac{1}{l^*}\sqrt{\frac{\Phi^*\epsilon^\sfE}{q C^*}}.
\end{displaymath}
We also introduce the following rescaled quantities,
\begin{equation}
\begin{array}{l}
t'=t/t^*,\ \ \bx'=\bx/l^*,\ \ \Phi'=\Phi/\Phi^*,\ \ \rho_z' =\rho_z/C*,\ \ z\in\{n,p,r,o,1,\cdots, N\}\\
R_{np}'(\rho_n',\rho_p')=\dfrac{t^*}{C^*}R_{np}(C^*\rho_n',C^*\rho_p'),\ \ \  R_{ro}'(\rho_r',\rho_o')=\dfrac{t^*}{C^*}R_{ro}(C^*\rho_r',C^*\rho_o')\\
R_{j}'(\rho_1',\cdots,\rho_N')=\dfrac{t^*}{C^*}R_{j}(C^*\rho_1',\cdots,C^*\rho_N'),\ \ \  1\le j\le N\\
G_0'=\dfrac{t^*}{l^* C^*}G_0,\ \ \ D_z'=D_z \dfrac{t^*}{{l^*}^2},\ \ \ \mu_z' =\mu_z\dfrac{t^*\Phi^*}{{l^*}^2},\ \ z\in\{n,p,r,o,1,\cdots, N\}\\
T'=T/t^*,\ \ C'=C/C^*,\ \ A_p'=t^*{C^*}^2 A_p',\ \ A_n'=t^*{C^*}^2 A_n,\ \ \rho_{isc}=\rho_{isc}/C^*
\end{array}
\end{equation}
for the model equations~\eqref{EQ:DDP} and~\eqref{EQ:DDP Redox}, and the following rescaled variables,
\begin{equation}
\begin{array}{llll}
k_{et}'=k_{et} t^* {C^*}/l^*, & k_{ht}'=k_{ht} t^* {C^*}/l^*, & v_n'=v_n t^*/l^*, & v_p'=v_p t^*/l^* \\
\varphi_{bi}'=\varphi_{bi}/\Phi^*, & \varphi_{app}'=\varphi_{app}/\Phi^*, & \varphi_{Stky}'=\varphi_{Stky}/\Phi^*, & \varphi_{app}^{\sfA'}=\varphi_{app}^\sfA/\Phi^*\\
\rho_n^{e'}=\rho_n^e/C^*, & \rho_p^{e'}=\rho_p^e/C^*, & \rho_z^{'\infty} =\rho_z^\infty/C^*, & z\in\{r,o,1,\cdots, N\}\\
\end{array}
\end{equation}
for the boundary and interface conditions described in~\eqref{EQ:DDP BC Dirich},~\eqref{EQ:DDP BC Robin},~\eqref{EQ:DDP BC Neumann},~\eqref{EQ:DDP Redox BC},~\eqref{EQ:DDP Redox BC Neumann},~\eqref{EQ:Interf E},~\eqref{EQ:Interf S},~\eqref{EQ:Interf Non Redox}, and~\eqref{EQ:Interf Phi}.

We can now summarize the mathematical model in rescaled (nondimensionalized) form as the following,
\begin{equation}\label{EQ:DDP ND}
\begin{array}{ll}
\partial_{t'} \rho_n' - \nabla\cdot (D_n'\nabla \rho_n' + \alpha_n \mu_n' \rho_n'\nabla\Phi') = -R_{np}'(\rho_n',\rho_p')+\gamma G_{np}', &\mbox{in}\ \  (0,T']\times\Omega_\sfS,\\
\partial_{t'} \rho_p' - \nabla\cdot (D_p'\nabla \rho_p' + \alpha_p\mu_p' \rho_p'\nabla\Phi') = -R_{np}'(\rho_n',\rho_p')+\gamma G_{np}', &\mbox{in}\ \  (0,T']\times\Omega_\sfS,\\
-\nabla\cdot{\lambda_\sfS^2}\nabla\Phi'= C' +\alpha_p\rho_p'+\alpha_n\rho_n',  &\mbox{in}\ \ (0,T']\times\Omega_\sfS\\
\Phi' = \varphi_{bi}'+\varphi_{app}' (\mbox{or}\ \ \Phi' = \varphi_{Stky}'+\varphi_{app}'),\ \ \rho_n'=\rho_n^{e'}, \ \ \rho_p'=\rho_p^{e'} & \mbox{on}\ (0,T']\times\Gamma_{\sfC}\\
\bnu\cdot\eps_r^\sfS\nabla\Phi' = 0,\ \  \bnu\cdot D_n'\nabla\rho_n'=0,\ \ \bnu\cdot D_p'\nabla\rho_p'=0 & \mbox{on}\ (0,T']\times\Gamma_{\sfS}\\
\Phi' = \varphi_{app}^{\sfA'}, & \mbox{on}\ (0,T']\times\Gamma_{\sfA}\\ 
\bnu\cdot\eps_r^\sfE\nabla\Phi' = 0, & \mbox{on}\ (0,T']\times\Gamma_{\sfE}
\end{array}
\end{equation}
for the transport dynamics in the rescaled semiconductor domain $\Omega_\sfS$, as the following,
\begin{equation}\label{EQ:DDP Redox ND}
\begin{array}{ll}
\partial_{t'} \rho_r' - \nabla\cdot (D_r'\nabla \rho_r' +\alpha_r \mu_r' \rho_r'\nabla\Phi') = +R_{ro}'(\rho_r',\rho_o'), &\mbox{in}\ \  (0,T']\times\Omega_\sfE,\\
\partial_{t'} \rho_o' - \nabla\cdot (D_o'\nabla \rho_o' +\alpha_o \mu_o' \rho_o'\nabla\Phi') = -R_{ro}'(\rho_r',\rho_o'), &\mbox{in}\ \  (0,T']\times\Omega_\sfE,\\
\partial_{t'} \rho_j' - \nabla\cdot (D_j'\nabla \rho_j' +\alpha_j \mu_j' \rho_j'\nabla\Phi') = R_{j}'(\rho_1',\cdots,\rho_N'), &\mbox{in}\ \  (0,T']\times\Omega_\sfE,\\
-\nabla\cdot{\lambda_\sfE^2}\nabla\Phi'= \alpha_o\rho_o'+\alpha_r\rho_r'+\sum_{j=1}^N \alpha_j\rho_j',  &\mbox{in}\ \ (0,T']\times\Omega_\sfE\\ 
\rho_r' = \rho_r^{'\infty},\ \ \rho_o' = \rho_o^{'\infty},\ \ \rho_j'=\rho_j^{'\infty}, & \mbox{on}\ \ (0,T']\times\Gamma_{\sfA}\\ 
\bnu\cdot D_r'\nabla \rho_r' =0,\ \ \bnu\cdot D_o'\nabla\rho_o' = 0,\ \ \bnu\cdot D_j'\nabla \rho_j'=0, & \mbox{on}\ \ (0,T']\times\Gamma_{\sfE}
\end{array}
\end{equation}
for the transport dynamics in the (rescaled) electrolyte domain $\Omega_\sfE$, and as the following,
\begin{equation}\label{EQ:DDP Interf ND}
\begin{array}{ll}
{[}\Phi'{]}_{\Sigma}=0,\ \  {[}\eps_r\frac{\partial\Phi'}{\partial\nu}{]}_{\Sigma}=0,\ \ \bnu\cdot D_j'\nabla \rho_j'=0 & \mbox{on}\ (0,T']\times \Sigma\\ 
\bnu^+\cdot\bJ_n'=-k_{et}'(\rho_n'-\rho_n^{e'})\rho_o',\ \  \bnu^+\cdot\bJ_p'=-k_{ht}'(\rho_p'-\rho_p^{e'})\rho_r', & \mbox{on}\ (0,T']\times\Sigma\\ 
\bnu\cdot\bJ_r'=-\bnu\cdot\bJ_o'=k_{ht}'(\rho_p'-\rho_p^{e'})\rho_r'-k_{et}'(\rho_n'-\rho_n^{e'})\rho_o', & \mbox{on}\ (0,T']\times\Sigma
\end{array}
\end{equation}
for the dynamics on the interface $\Sigma$. Here the rescaled Auger generation-recombination rate $R_{np}'$ takes exactly the same form as in~\eqref{EQ:Auger}, i.e.,
\begin{equation}\label{EQ:Auger ND}
	R_{np}'(\rho_n',\rho_p') = (A_n' \rho_n' +A_p' \rho_p')(\rho_{isc}^{'2}-\rho_n'\rho_p'),
\end{equation}
and the rescaled photon illumination function $G_{np}'$ takes the form
\begin{equation}\label{EQ:Illu Source ND}
	G_{np}'(\bx')=\left\{
	\begin{array}{rr}
		\sigma'(\bx') G_0'(\bx_0')\exp\Big(-\int_0^{\bar s}\sigma'(\bx_0'+s\btheta_0)ds\Big),& \mbox{if}\ \bx'=\bx_0'+\bar s\btheta_0\\
	0,& \mbox{otherwise}
	\end{array}\right.
\end{equation}
In the rest of the paper, we will work on the numerical simulations of the semiconductor-electrolyte system based on the above nondimensionalized systems~\eqref{EQ:DDP ND} and~\eqref{EQ:DDP Redox ND}.

\subsection{Time-dependent discretization}

For the numerical simulations in this paper, we discretize the time-dependent systems of reaction-drift-diffusion-Poisson equations~\eqref{EQ:DDP ND} and~\eqref{EQ:DDP Redox ND} by standard finite difference method in both spatial and temporal variables. In the spatial variable, we employ a classical upwind discretization of the advection terms (such as $\nabla\Phi'\cdot\nabla\rho_n'$) to ensure the stability of the scheme. To avoid solving nonlinear systems of equations in each time step (since the models~\eqref{EQ:DDP ND} and~\eqref{EQ:DDP Redox ND} are nonlinear), we employ the forward Euler scheme for the temporal variable. Since this is a first-order scheme and is explicit, we do not need to perform any nonlinear solve in the solution process, as long as we can supply the right initial conditions. We are aware that there are many efficient solvers for similar problems that have been developed in the literature; see, for instance, ~\cite{ZhOvLiZhNi-JCP11}. 

\medskip

To solve stationary problems, we can evolve the system for a long time so that the system reaches its stationary state. We use the magnitude of the relative $L^2$ update of the solution as the stopping criterion. An alternative, in fact more efficient, way to solve the nonlinear system is the following iterative scheme.

\subsection{A Gummel-Schwarz iteration for stationary problems}

This method combines domain decomposition strategies with nonlinear iterative schemes. We decompose the system naturally into two subsystems, the semiconductor system and the electrolyte system. We solve the two subsystem alternatively and couple them with the interface condition. This is the Schwarz decomposition strategy that has been used extensively in the literature; see~\cite{ChELiSh-JCP07,MiQuSa-JCP95} for similar domain decomposition strategies in semiconductor simulation. To solve the nonlinear equations in each sub-problem, we adopt the Gummel iteration scheme~\cite{BeCaCaVe-JCP09,BuPi-M3AS09,Kulikovsky-JCP95}. This scheme decomposes the drift-diffusion-Poisson system into a drift-diffusion part and a Poisson part and then solves the two parts alternatively. The coupling then comes from the source term in the Poisson equation. Our algorithm, in the form of solving the stationary problem, takes the following form.\\

\noindent{\bf GUMMEL-SCHWARZ ALGORITHM.}
\begin{itemize}
\item[{[1]}] Gummel step $k=0$, construct initial guess $\{\rho_n^{'0}, \rho_p^{'0}, \rho_r^{'0}, \rho_o^{'0}, \{\rho_j^{'0}\}_{j=1}^N\}$
\item[{[2]}] Gummel step $k\ge 1$:
	\begin{itemize}
	\item Solve the Poisson problem for $\Phi^{'k}$ in $\Omega_\sfE\cup\Omega_\sfS$ using the densities $\rho_n^{'k-1}, \rho_p^{'k-1}, \rho_r^{'k-1}$, $\rho_o^{'k-1}$, and $\{\rho_j^{'k-1}\}_{j=1}^N$:
	\begin{equation}\label{EQ:DDP ND Poisson}
	\begin{array}{ll}
		-\nabla\cdot{\lambda_\sfS^2}\nabla\Phi^{'k}= C' +\alpha_p\rho_p^{'k-1}+\alpha_n\rho_n^{'k-1},  &\text{in}\ \Omega_\sfS \\
		-\nabla\cdot{\lambda_\sfE^2}\nabla\Phi^{'k}= \alpha_o\rho_o^{'k-1}+\alpha_r\rho_r^{'k-1}+\sum_{j=1}^N\alpha_j\rho_j^{'k-1},  &\mbox{in}\ \Omega_\sfE\\
		\Phi^{'k} = \varphi_{bi}'+\varphi_{app}',\ \ \mbox{or}\ \ \Phi^{'k} = \varphi_{Stky}'+\varphi_{app}', & \mbox{on}\ \Gamma_{\sfC}\\
		\bnu\cdot\eps_r^\sfS\nabla\Phi^{'k} = 0, & \mbox{on}\ \Gamma_{\sfS}\\
		\Phi^{'k} = \varphi_{app}^{\sfA'}, & \mbox{on}\ \Gamma_{\sfA}\\
		\bnu\cdot\eps_r^\sfE\nabla\Phi^{'k} = 0, & \mbox{on}\ \Gamma_{\sfE}\\
	{[}\Phi^{'k}{]}_{\Sigma}=0,\ \  {[}\eps_r\frac{\partial\Phi^{'k}}{\partial\nu}{]}_{\Sigma}=0,& \mbox{on}\ \Sigma
	\end{array}
	\end{equation}
	\item Solve for $\rho_n^{'k},\rho_p^{'k}$, $\rho_r^{'k},\rho_o^{'k},\{\rho_j^{'k}\}_{j=1}^N$ as limits of the following iteration:
	\begin{itemize} 
		\item[{[i]}] Schwarz step $\ell=0$: construct guess $(\rho_n^{'k,0},\rho_p^{'k,0},\rho_r^{'k,0},\rho_o^{'k,0},\{\rho_j^{'k,0}\}_1^N)$
		\item[{[ii]}] Schwarz step $\ell\ge 1$: solve \emph{sequentially}
		\begin{equation}\label{EQ:SchI}
		\begin{array}{ll}
			\nabla\cdot (-D_n'\nabla \rho_n^{'k,\ell} + \alpha_n\mu_n' \rho_n^{'k,\ell}\nabla\Phi^{'k}) = -R_{np}'(\rho_n^{'k,\ell},\rho_p^{'k,\ell})+\gamma G_{np}', &\mbox{in}\ \Omega_\sfS\\
\nabla\cdot (-D_p'\nabla \rho_p^{'k,\ell} +\alpha_p \mu_p' \rho_p^{'k,\ell}\nabla\Phi^{'k})= -R_{np}'(\rho_n^{'k,\ell},\rho_p^{'k,\ell})+\gamma G_{np}', &\mbox{in}\ \Omega_\sfS\\ \rho_n^{'k,\ell}=\rho_n^{e'}, \ \ \rho_p^{'k,\ell}=\rho_p^{'e}, & \mbox{on}\ \Gamma_{\sfC}\\ 
\bnu\cdot D_n'\nabla\rho_n^{'k,\ell}=0,\ \ \bnu\cdot D_p'\nabla\rho_p^{'k,\ell}=0, & \mbox{on}\ \Gamma_{\sfS}\\ 
\bnu^+\cdot\bJ_n^{'k,\ell}=-k_{et}'(\rho_n^{'k,\ell}-\rho_n^{e'})\rho_o^{'k,\ell-1},\ \  \bnu^+\cdot\bJ_p^{'k,\ell}=-k_{ht}'(\rho_p^{'k,\ell}-\rho_p^{e'})\rho_r^{'k,\ell-1}, & \mbox{on}\ \Sigma
\end{array}
\end{equation}
and
		\begin{equation}\label{EQ:SchII}
		\begin{array}{ll}
			\nabla\cdot(-D_r'\nabla \rho_r^{'k,\ell} + \alpha_r\mu_r' \rho_r^{'k,\ell}\nabla\Phi^{'k}) 					= +R_{ro}'(\rho_r^{'k,\ell},\rho_o^{'k,\ell}), &\mbox{in}\ \Omega_\sfE\\
			\nabla\cdot(-D_o'\nabla \rho_o^{'k,\ell} +\alpha_o \mu_o' \rho_o^{'k,\ell}\nabla\Phi^{'k}) 					= -R_{ro}'(\rho_r^{'k,\ell},\rho_o^{'k,\ell}), &\mbox{in}\ \Omega_\sfE\\
			\nabla\cdot(-D_j'\nabla \rho_j^{'k,\ell} +\alpha_j \mu_j' \rho_j^{'k,\ell}\nabla\Phi^{'k}) 					= R_{j}(\rho_1^{k,\ell},\cdots,\rho_N^{'k,\ell}), &\mbox{in}\ \Omega_\sfE\\
		\rho_r^{'k,\ell} = \rho_r^{'\infty},\ \ \rho_o^{'k,\ell} = \rho_o^{'\infty},\ \ \rho_j^{'k,\ell}=\rho_j^{'\infty}, & \mbox{on}\ \Gamma_{\sfA}\\
		\bnu\cdot D_r'\nabla \rho_r^{'k,\ell} =0,\ \ \bnu\cdot D_o'\nabla\rho_o^{'k,\ell} = 0,\ \ \bnu\cdot D_j'\nabla \rho_j^{'k,\ell}=0, & \mbox{on}\ \ \Gamma_{\sfE}\\
		\bnu\cdot\bJ_r^{'k,\ell}=-\bnu\cdot\bJ_o^{'k,\ell}=k_{ht}'(\rho_p^{'k,\ell}-\rho_p^{e'})\rho_r^{'k,\ell}-k_{et}'(\rho_n^{'k,\ell}-\rho_n^{e'})\rho_o^{'k,\ell}, & \mbox{on}\ \Sigma\\
		\bnu\cdot D_j'\nabla \rho_j^{'k,\ell}=0,\ \ 1\le j\le N & \mbox{on}\ \ \Sigma.
		\end{array}
		\end{equation}
		\item[{[iii]}] If convergence criteria satisfied, stop; otherwise, set $\ell=\ell+1$ and go to [ii].
	\end{itemize}
	\end{itemize}
\item[{[3]}] If convergence criteria satisfied, stop; Otherwise, set $k=k+1$ and go to step [2].
\end{itemize}
Note that since~\eqref{EQ:SchI} and ~\eqref{EQ:SchII} are solved \emph{sequentially}, we are able to replace the $\rho_n^{'k,\ell-1}$ and $\rho_p^{'k,\ell-1}$ terms in~\eqref{EQ:SchII} in the boundary conditions on $\Sigma$ with $\rho_n^{'k,\ell}$ and $\rho_p^{'k,\ell}$ respectively. This slightly improves the speed of convergence of the numerical scheme.
If the mathematical system is well-posed, the convergence of this iteration can be established following the lines of work in~\cite{Lions-DDM87,Lions-DDM89}. The details will be in a future work. 

\subsection{Computational issues}

\paragraph{Initial conditions.} To solve the time-dependent problem, we need to supply the system with appropriate initial conditions. Our initial conditions are given only for the densities. The initial condition on the potential is obtained by solving the Poisson equation with appropriate boundary conditions using the initial densities. That is done by solving problem~\eqref{EQ:DDP ND Poisson}.

\paragraph{Stiffness of the system.} One major issue in the numerical solution of the mathematical model is the stiffness of the system across the interface, caused by the large contrast between the magnitudes of the densities in the semiconductor and electrolyte domains. Physically this leads to the formation of boundary layers of charge densities local to the interface; see, for instance, discussions in~\cite{BaFa-Book00,Memming-Book01}. To capture the sharp transition at the interface, we have to use finer finite difference grids around the interface. 

\paragraph{Computation of stationary states.} The Gummel-Schwarz iteration is attractive when stationary state solutions are to be sought since it is computationally much more efficient than the time-stepping method. However, due to the nonlinearity of the system, multiple solutions could exist, as illustrated, for instance, in~\cite{WaReOd-SIAM91} in similar settings. Starting from different initial guesses for the nonlinear solver could lead us to different solutions (and some of them could be unphysical). The time-stepping scheme, however, will lead us to the desired steady state for any given initial state. We can combine the two methods. For a given initial state, we use the time-stepping scheme to evolve the system for some time, say $\tilde T'$ (which we select empirically as a time when the fastest initial evolution of the solution has passed); we then use the solution at $\tilde T'$ as the initial guess for the Gummel-Schwarz iteration to find the corresponding steady state solution. In our numerical simulations, we did not observe the existence of multiple steady states. However, this is an issue to keep in mind when more complicated (for instance multidimensional) simulations are performed.

\section{Numerical simulations}
\label{SEC:Num}

We present some numerical simulations based on the mathematical model that we have constructed for the semiconductor-electrolyte system.

\subsection{General setup for the simulations}

To simplify the numerical computation, we assume some symmetry in the semiconductor-electrolyte system so that we can reduce the problem to one dimension. We showed in Fig.~\ref{FIG:Semi-Liq Cell} (middle and right) two typical two-dimensional systems where such dimension reductions can be performed. In the first case, if we assume that the system is invariant in the direction parallel to the current collector, then we have a one-dimensional system in the direction that is perpendicular to the current collector. In the second setting, we have a radially symmetric system that is invariant in the angular direction in the polar coordinate. The system can then be regarded as a one-dimensional system in the radial direction. 

In all the simulations we have performed, we use an $n$-type semiconductor to construct the semiconductor-electrolyte system. We also select the electrolyte such that the charge numbers $\alpha_o-\alpha_r=1$. Furthermore, we assume in our simulation that the reductants and oxidants are the only charges in the electrolyte. Therefore the equations for $\rho_j'$ $(1\le j\le N)$ are dropped. In Tab.~\ref{TAB:Para}, we list the values of all the model parameters that we use in the simulations. They are mainly taken from ~\cite{Brennan-Book99,KiBe-Book91,Schenk-Book98,YuCa-Book03} and references therein. These values may depend on the materials used to construct the system and the temperature of the system (which we set to be $\cT=300K$). They can be tuned to be more realistic by careful calibrations.
\begin{table}[ht]
  \centering
        \begin{tabular}{lll|lll}\hline
        Parameter& Value & Unit & Parameter & Value & Unit\\ \hline
        $q$          & $1.6 \times 10^{-19}$     & [${\rm A\ s}$] & 
	$k_B$        & $8.62\times 10^{-5}\ q$   & [${\rm J\ K^{-1}}$]\\ \hline
        $\eps_0$     & $8.85\times 10^{-14}$     & [${\rm A\ s\ V^{-1}\ cm^{-1}}$] & 
	$\eps_r^\sfS$     & $11.9$                    & \\ \hline
        $\mu_n$      & $1500$                    & [${\rm cm^2\ V^{-1}\ s^{-1}}$] & 
	$\mu_p$      & $450$                     & [${\rm cm^2\ V^{-1}\ s^{-1}}$]\\ \hline
        $A_n$    & $2.8\times 10^{-31}$      & [${\rm cm^6\ s^{-1}}$] & 
	$A_p$    & $9.9\times 10^{-32}$      & [${\rm cm^6\ s^{-1}}$]\\ \hline
	$N_c$    & $2.80\times 10^{19}$      & [${\rm cm^{-3}}$] & 
	$N_v$    & $1.04\times 10^{19}$      & [${\rm cm^{-3}}$]\\ \hline
	$v_n$    & $5\times 10^6$            & [${\rm cm\ s^{-1}}$] & 
	$v_p$        & $5\times 10^6$            & [${\rm cm\ s^{-1}}$]\\ \hline
	$\Phi_m$ & $2.4$  & [${\rm V}$] & $\chi$  & $1.2$ & [${\rm V}$]\\ \hline
	$k_{et}$        & $1\times 10^{-21}$        & [${\rm cm^4\ s^{-1}}$] &
	$k_{ht}$        & $1\times 10^{-17}$        & [${\rm cm^4\ s^{-1}}$]\\ \hline
	$\mu_o$      & $2\times 10^{-1}$                    & [${\rm cm^2\ V^{-1}\ s^{-1}}$]&
	$\mu_r$      & $0.5\times 10^{-1}$                       & [${\rm cm^2\ V^{-1}\ s^{-1}}$]\\ \hline
	$\eps_r^\sfE$ & $1000$                    & & $G_0$ & $1.2\times 10^{17}$ & [${\rm cm}^{-2} {\rm s}^{-1}$]\\ \hline
        \end{tabular}
        \caption{Values of physical parameters used in the numerical simulations. The numbers are given in unit of {\rm cm}, {\rm s}, {\rm V}, {\rm A}.}
\label{TAB:Para}
\end{table}

We perform simulations on two devices of different sizes. The two devices are designed to mimic a large (Device I) and small (Device II) nano- to micro-scale solar cell building block, as follows.

\medskip

\noindent{\bf Device I.} The device is contained in $\Omega=(-1.0, 1.0)$ with the semiconductor $\Omega_{\sfS}=(-1.0, 0)$ and $\Omega_{\sfE}=(0, 1.0)$ separated by the interface $\Sigma$ located at $x'=0$. The semiconductor boundary $\Gamma_{\sfS}$ is thus at the point $x'=-1.0$, while the electrolyte boundary $\Gamma_{\sfE}$ is at the point $x'=1.0$.

\medskip
\noindent{\bf Device II.} The device is contained in $\Omega=(-0.2,0.2)$ with the semiconductor $\Omega_{\sfS}=(-0.2, 0)$ and $\Omega_{\sfE}=(0, 0.2)$ separated by the interface $\Sigma$ located at $x'=0$. The semiconductor boundary $\Gamma_{\sfS}$ is thus at the point $x'=-0.2$, while the electrolyte boundary $\Gamma_{\sfE}$ is at the point $x'=0.2$.

\subsection{General dynamics of semiconductor-electrolyte systems}

We now present simulation results on general dynamics of the semiconductor-electrolyte systems we have constructed. We perform the simulations using the first system, i.e., Device I. General parameters in the simulations are listed in Tab.~\ref{TAB:Para}. We consider three different cases.

\begin{figure}[ht]
\centering
\includegraphics[height=0.3\textheight,width=0.45\textwidth]{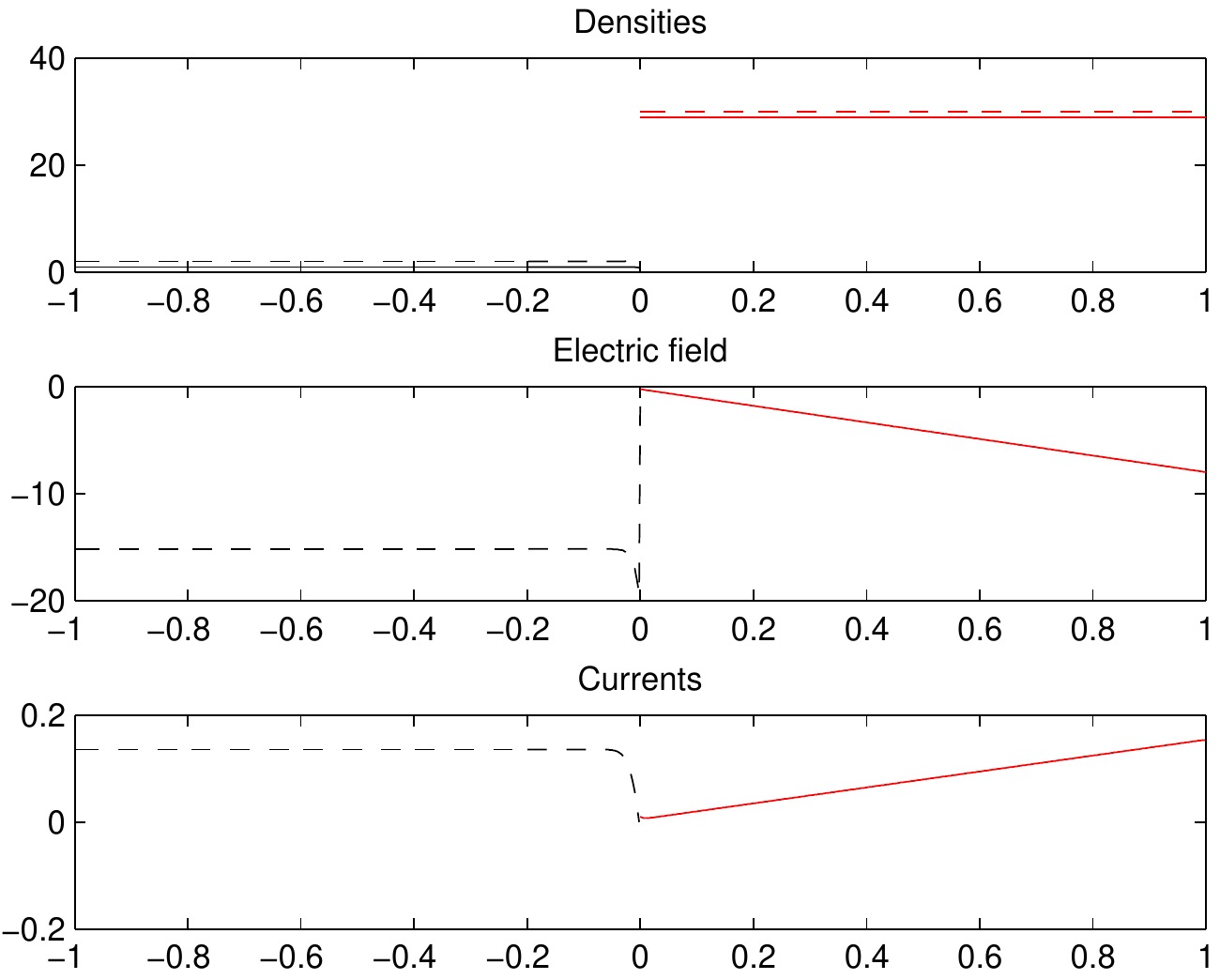}	\includegraphics[height=0.3\textheight,width=0.45\textwidth]{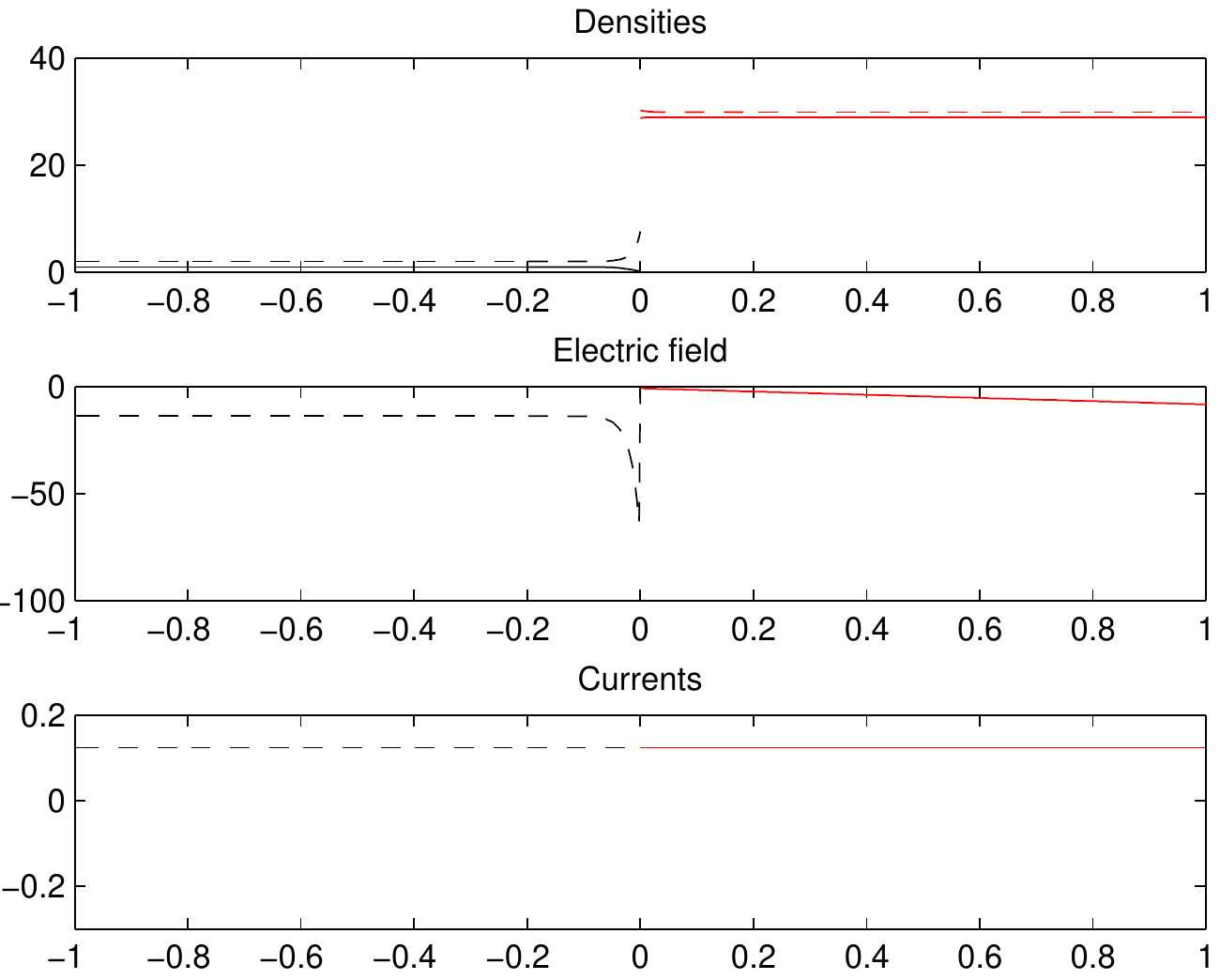}
\caption{Case I (a) in dark  environment ($\gamma=0$). Left column: charge densities (black solid for $\rho_p'$, black dashed for $\rho_n'$, red solid for $\rho_o'$, and red dashed for $\rho_r'$), electric field, and flux distributions at time $t'=0.05$. Right column: charge densities, electric field, and flux distributions at the stationary state.}
\label{FIG:Dark I(a)}
\end{figure}
\paragraph{Case I (a).} In this simulation, we show typical dynamics of the semiconductor-electrolyte system in dark and illuminated cases. We consider the case when the densities of the reductant-oxidant pair are very high compared to the densities of the electron-hole pair. Specifically, we take $\rho_r^{'\infty}= 30.0$, $\rho_o^{'\infty}=29.0$. The system starts from the initial conditions: $\rho_n^{'0}=2$, $\rho_p^{'0}=1.0$, $\rho_r^0 = \rho_r^{'\infty}$, and $\rho_o^{'0}  = \rho_o^{'\infty}$. We first perform simulation for the dark case, i.e., when the parameter $\gamma=0$. In the left column of Fig.~\ref{FIG:Dark I(a)}, we show the distributions of the charge densities, electric field and total fluxes in the semiconductor and the electrolyte at time $t'=0.05$. The corresponding distributions at stationary state are shown in the right column. 

\begin{figure}[ht]
\centering
\includegraphics[height=0.3\textheight,width=0.45\textwidth]{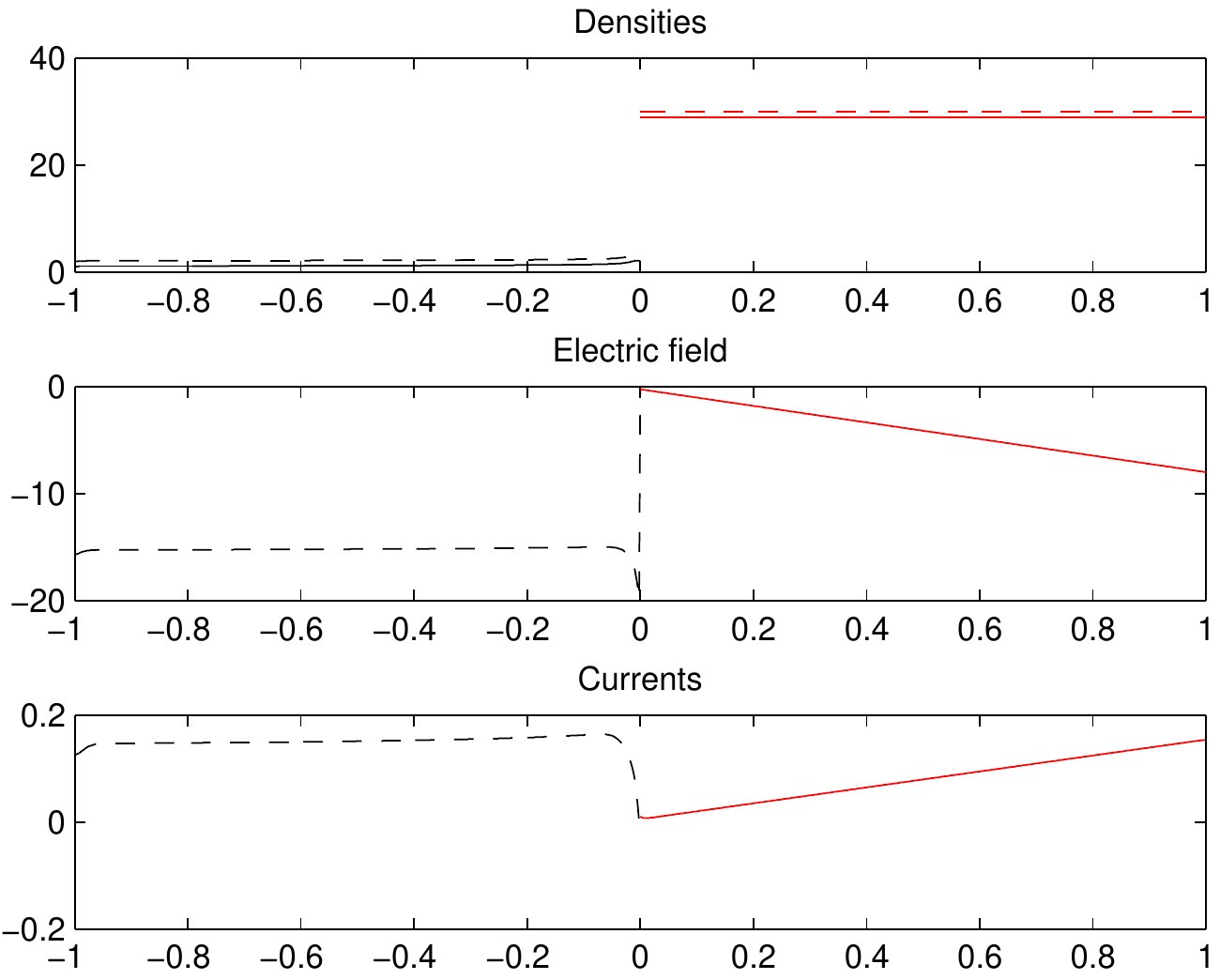} 	\includegraphics[height=0.3\textheight,width=0.45\textwidth]{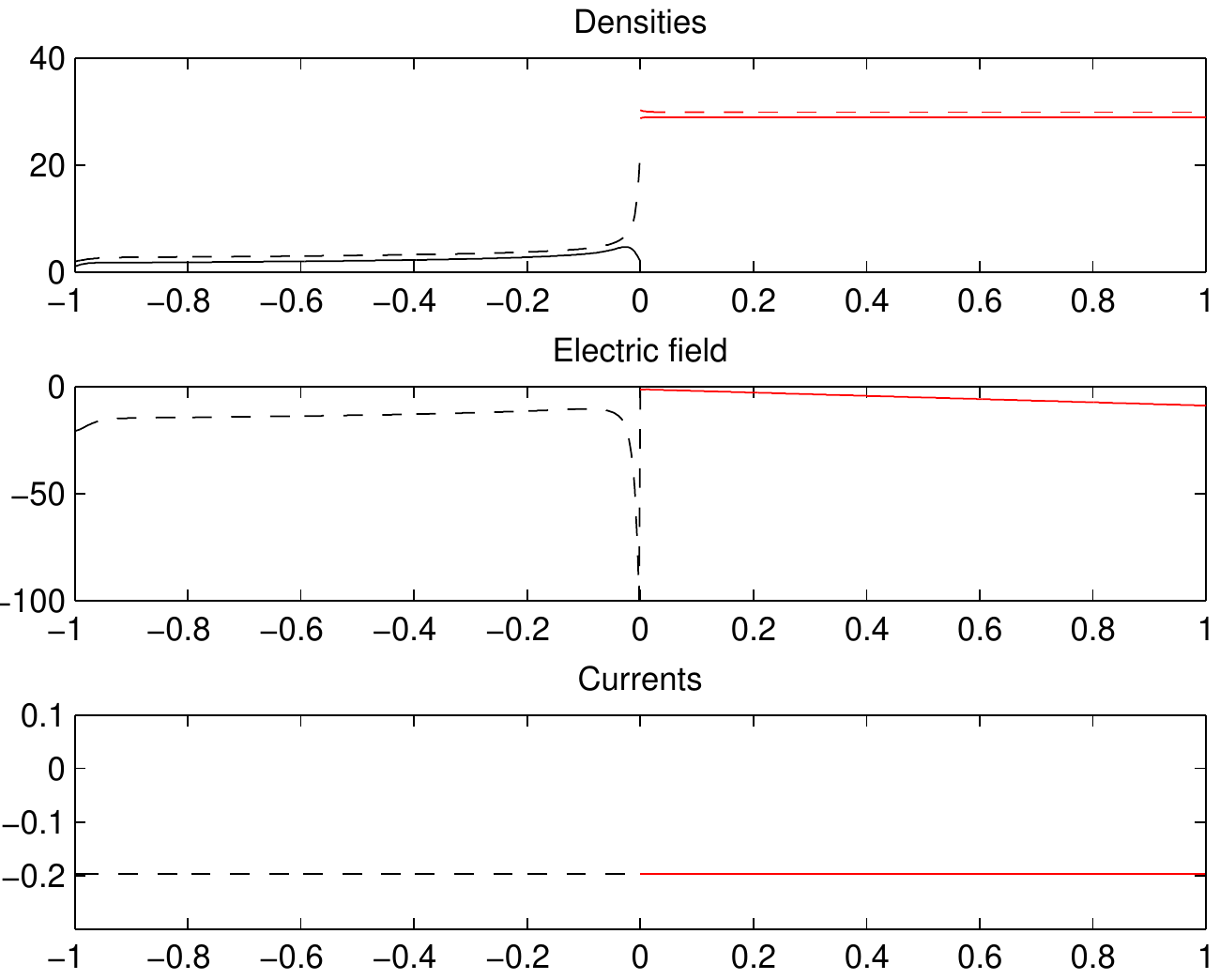}
\caption{Case I (a) in illuminated environment ($\gamma=1$). Left column: charge densities (black solid for $\rho_p'$, black dashed for $\rho_n'$, red solid for $\rho_o'$, and red dashed for $\rho_r'$), electric field, and flux distributions at time $t'=0.05$. Right column: charge densities, electric field and flux distributions at the stationary state.}
\label{FIG:Ill I(a)}
\end{figure}
We now repeat the simulation in the illuminated environment by setting the parameter $\gamma=1$. In this case, the surface photon flux is set as $G_0'=1.2\times 10^{-7}$ to mimic a solar spectral irradiance with air mass 1.5. The simulation results are presented in Fig.~\ref{FIG:Ill I(a)}. In both the dark and illuminated cases shown here, the applied potential bias is $19.3$, and Ohmic contact (i.e. Dirichlet condition) is assumed at the left end of the semiconductor.

\begin{figure}[hbt]
\centering 
\includegraphics[height=0.3\textheight,width=0.45\textwidth]{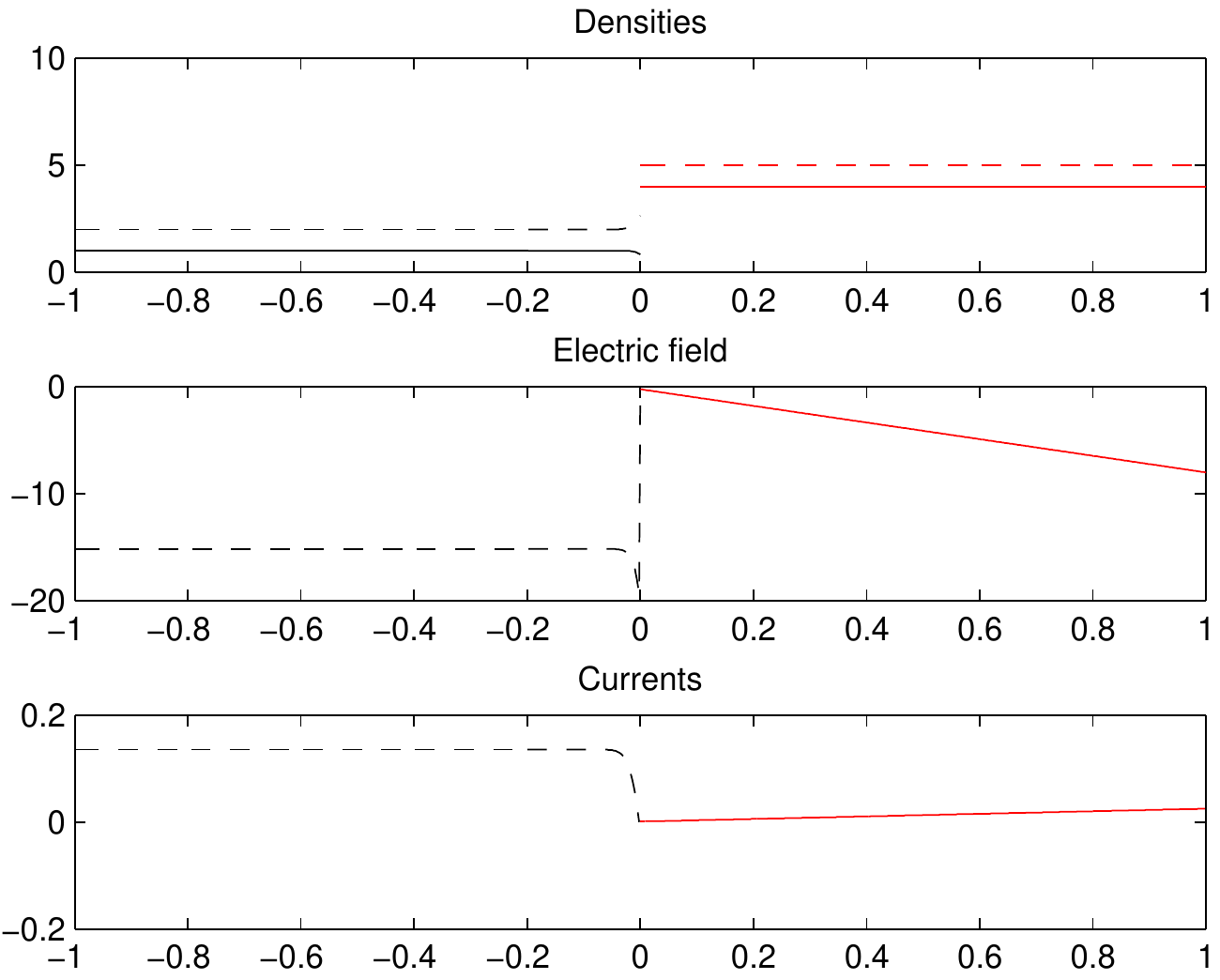}
\includegraphics[height=0.3\textheight,width=0.45\textwidth]{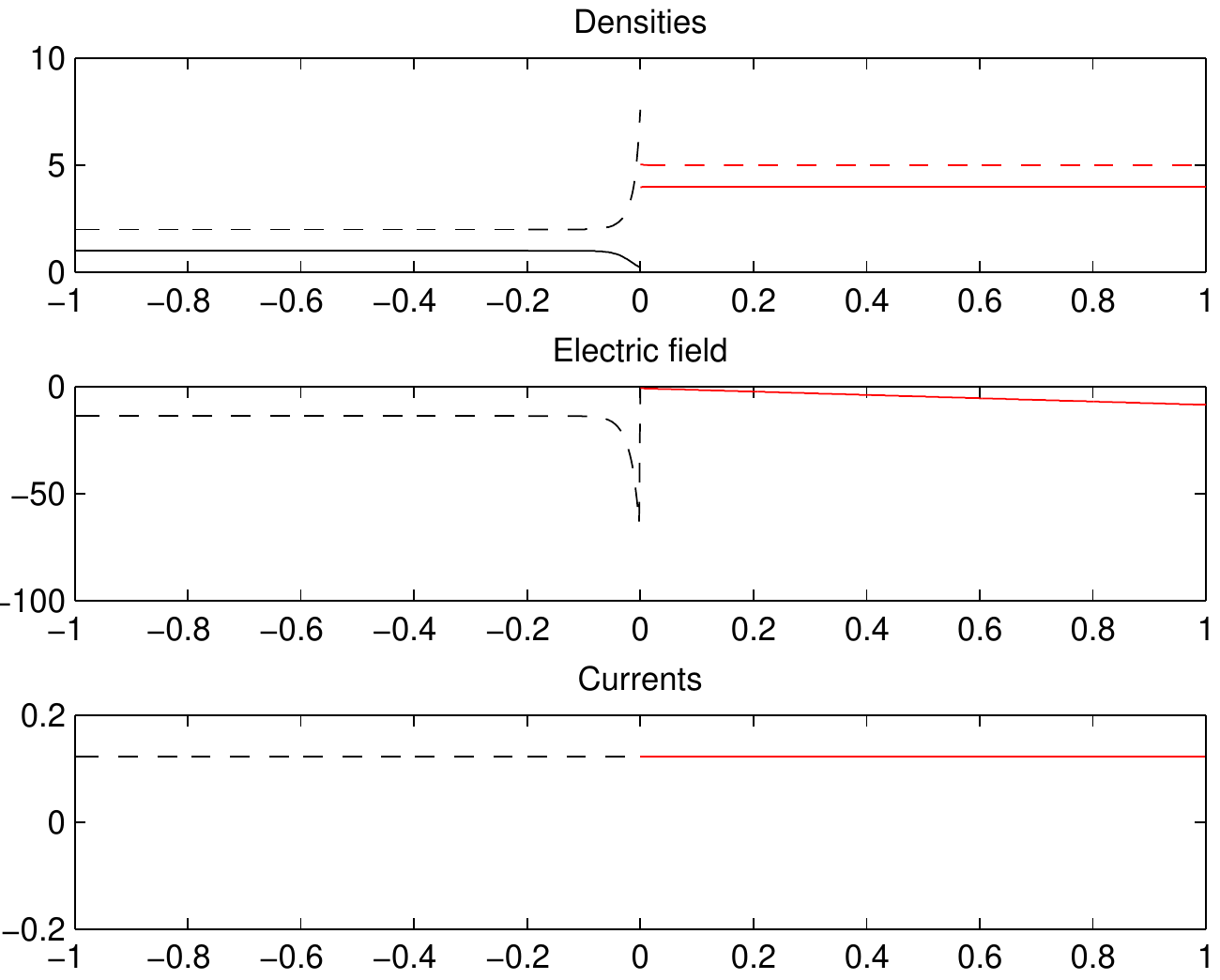}
\caption{Case I (b) in dark environment ($\gamma=0$). Left column: charge densities (black solid for $\rho_p'$, black dashed for $\rho_n'$, red solid for $\rho_o'$, and red dashed for $\rho_r'$), electric field, and flux distributions at time $t'=0.05$. Right column: charge densities, electric field and flux distributions at the stationary state.}
\label{FIG:Dark I(b)}
\end{figure}
\begin{figure}[hbt]
\centering
\includegraphics[height=0.3\textheight,width=0.45\textwidth]{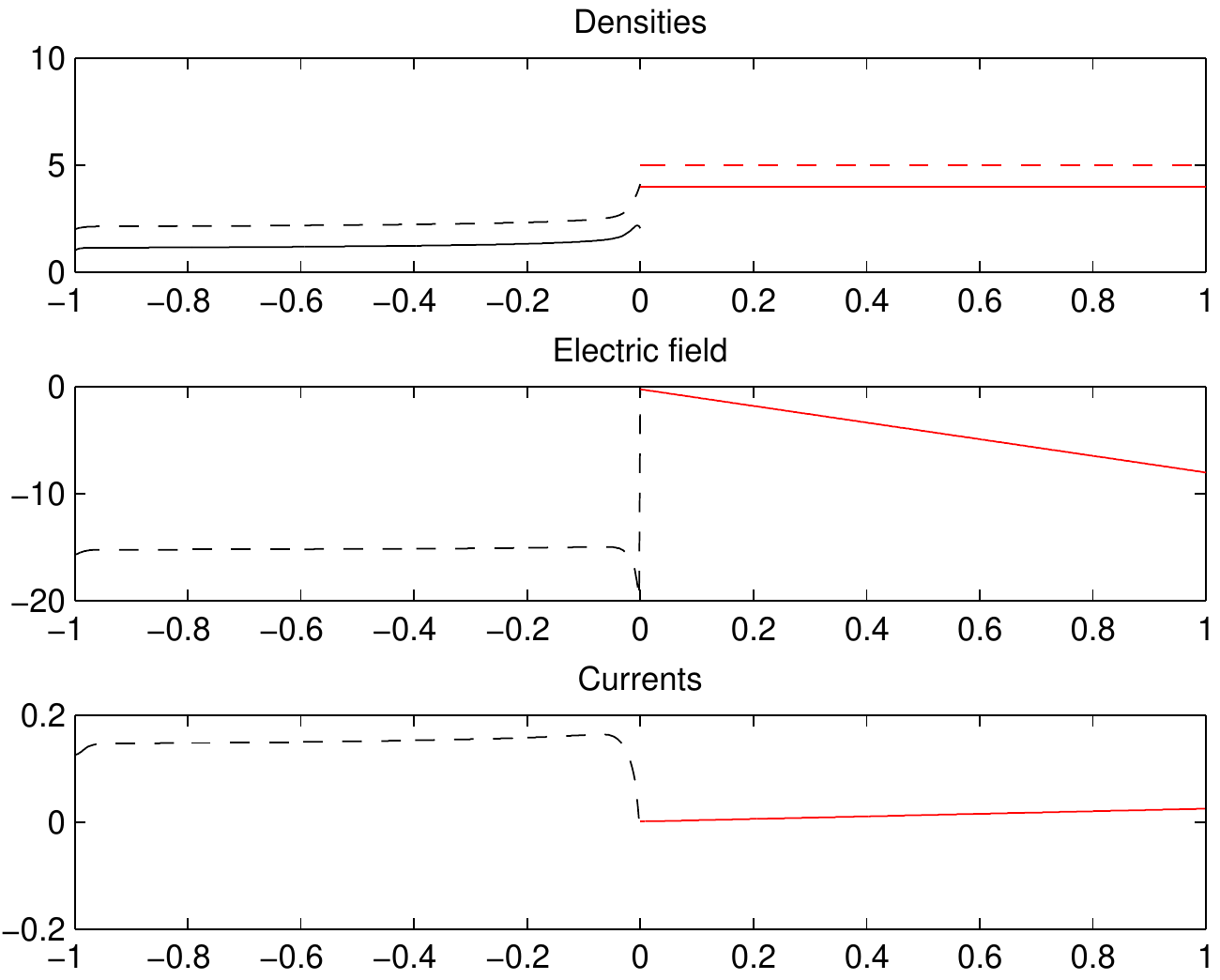}
\includegraphics[height=0.3\textheight,width=0.45\textwidth]{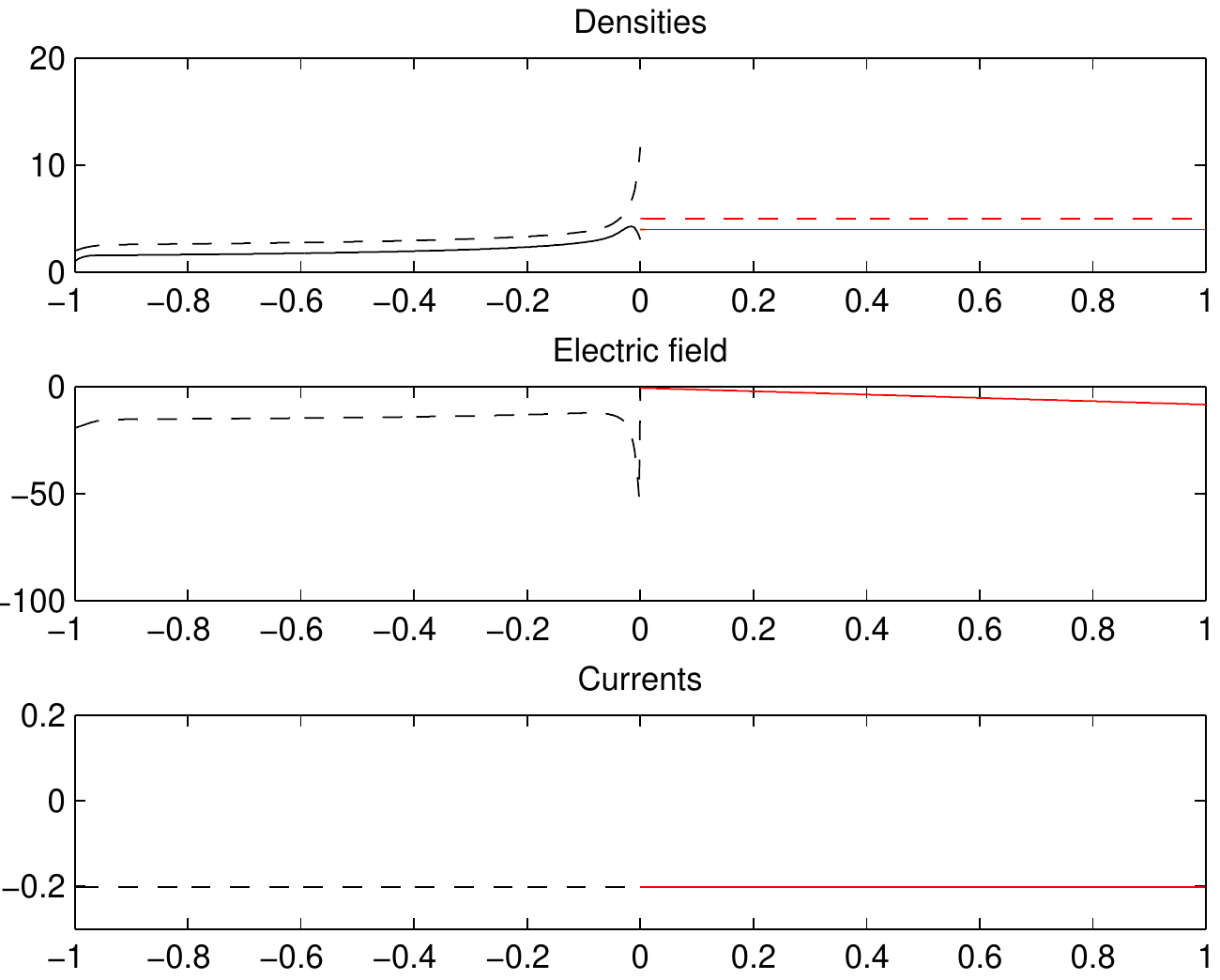}
\caption{Case I (b) in illuminated environment ($\gamma=1$). Left column: charge densities (black solid for $\rho_p'$, black dashed for $\rho_n'$, red solid for $\rho_o'$, and red dashed for $\rho_r'$), electric field, and flux distributions at time $t'=0.05$. Right column: charge densities, electric field, and flux distributions at the stationary state.}
\label{FIG:Ill I(b)}
\end{figure}

\paragraph{Case I (b).} We now perform a set of simulations with the same parameters as in the previous numerical experiment but with a lower bulk reductant and oxidant pair densities: $\rho_r^{'\infty}  = 4.0$, $\rho_o^{'\infty}= 5.0$. The initial conditions for the simulation are: $\rho_n^{'0}=2.5$, $\rho_p^{'0}=1.0$, $\rho_r^{'0} = \rho_r^{'\infty}$ and $\rho_o^{'0}  = \rho_o^{'\infty}$. The results in the dark and the illuminated environments are shown in Fig.~\ref{FIG:Dark I(b)} and Fig.~\ref{FIG:Ill I(b)}, respectively. 
We observe from the comparison of Fig.~\ref{FIG:Dark I(a)} and Fig~\ref{FIG:Ill I(a)} with Fig.~\ref{FIG:Dark I(b)} and Fig.~\ref{FIG:Ill I(b)}, that when all other factors are kept unchanged, lowering the density of the redox pair leads to significant change of the electric field across the device, especially at the semiconductor-electrolyte interface. The charge densities in the semiconductor also change significantly. In both cases, however, the charge densities in the electrolyte remain as almost constant across the electrolyte. There are two reasons for this. First, the charge densities in the electrolyte are significantly higher than those in the semiconductor (even in Case I (b)). Second, the relative dielectric constant of the electrolyte is much higher than that in the semiconductor, resulting in a relatively constant electric field inside the electrolyte. If we lower the relative dielectric function, we observe a significant variation in charge density distributions in the electrolyte, as we can see from the next numerical simulation.

\begin{figure}[htb]
\centering 
\includegraphics[height=0.3\textheight,width=0.45\textwidth]{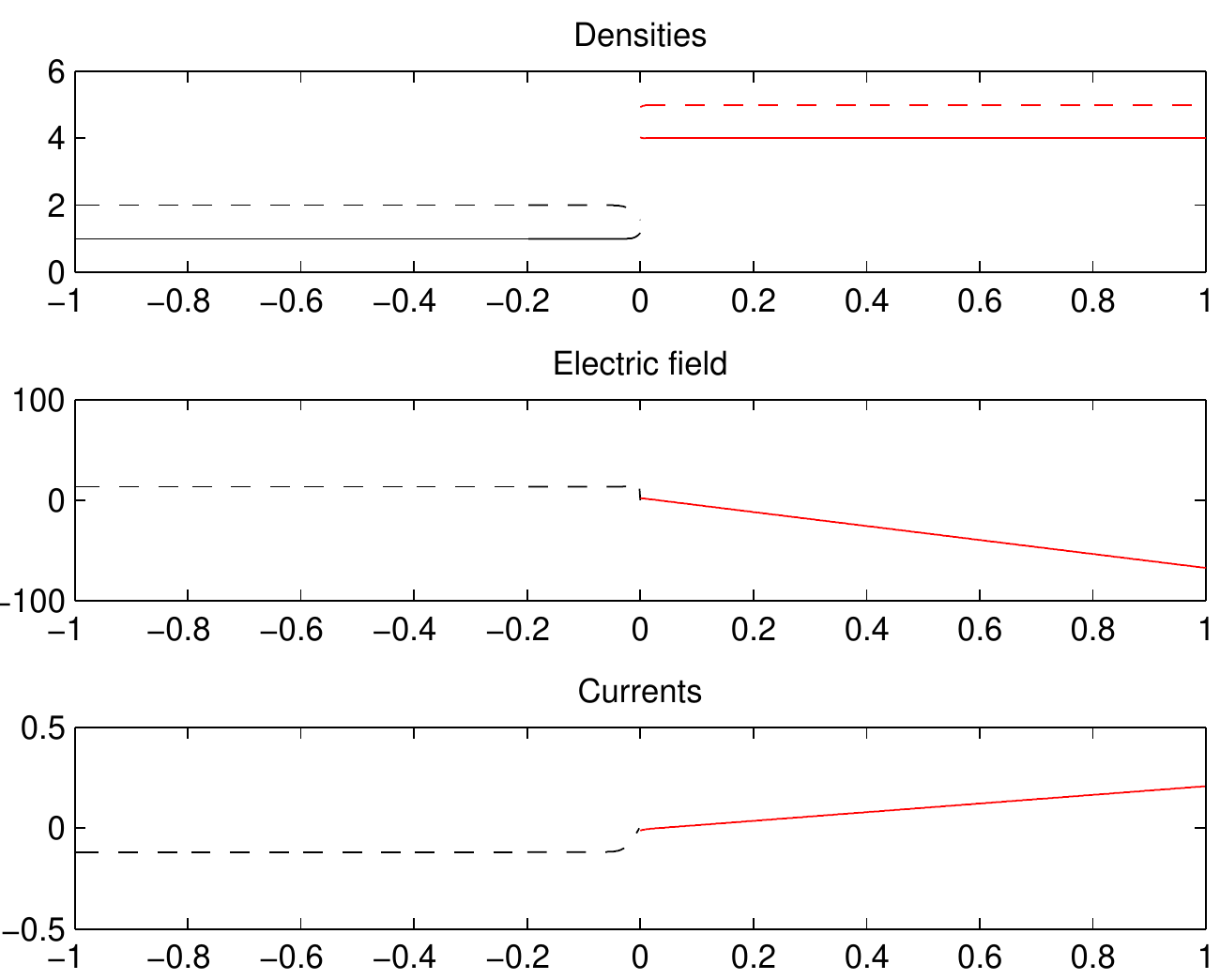} 
\includegraphics[height=0.3\textheight,width=0.45\textwidth]{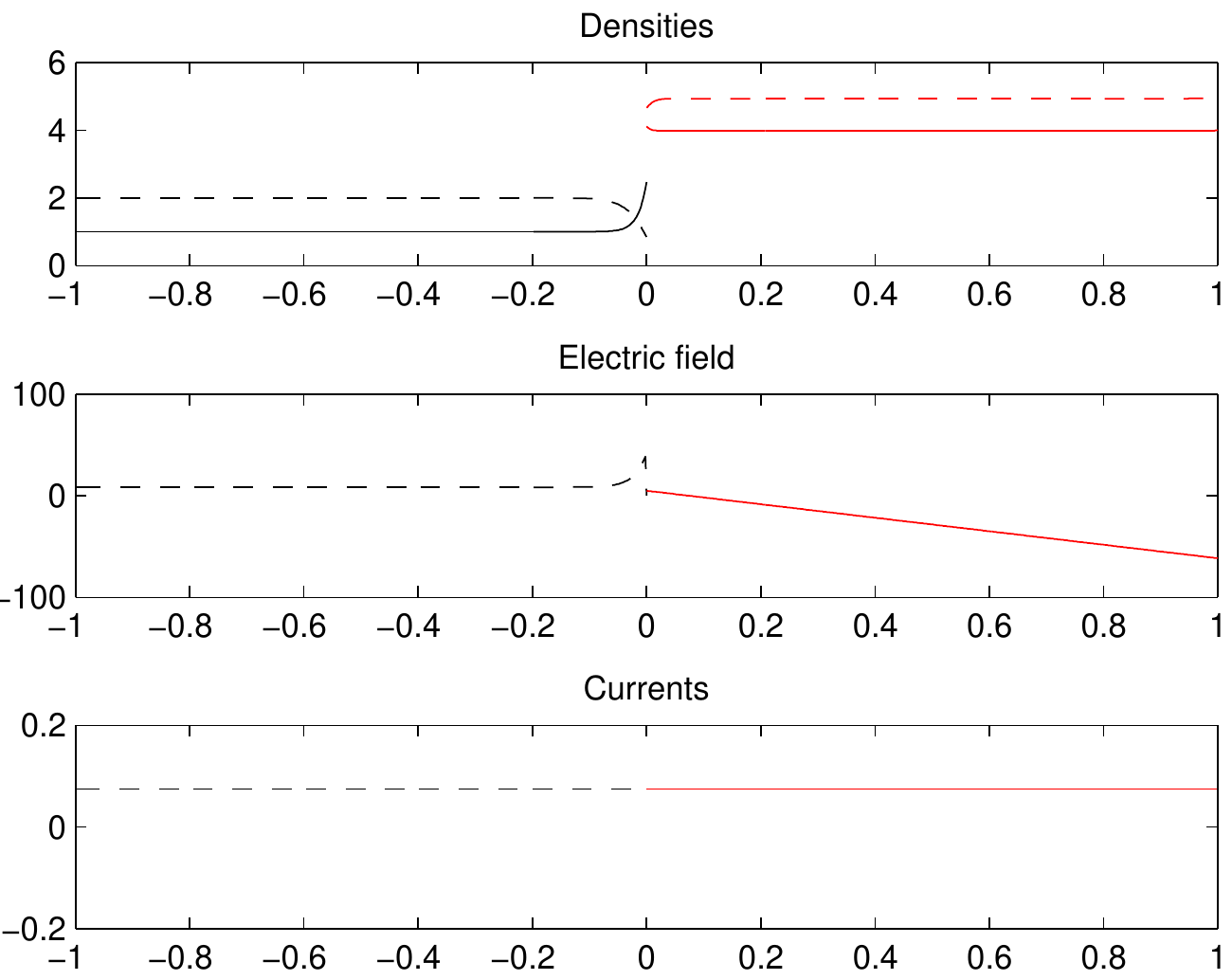}
\caption{Case I (c) in dark environment ($\gamma=0$). Shown are: charge densities (black solid for $\rho_p'$, black dashed for $\rho_n'$, red solid for $\rho_o'$, and red dashed for $\rho_r'$), electric field, and flux distributions at time $t'= 0.05$ (left); and charge densities, electric field, and flux distributions at the stationary state (right).}
\label{FIG:Dark I(c)}
\end{figure}
\begin{figure}[htb]
\centering 
\includegraphics[height=0.3\textheight,width=0.4\textwidth]{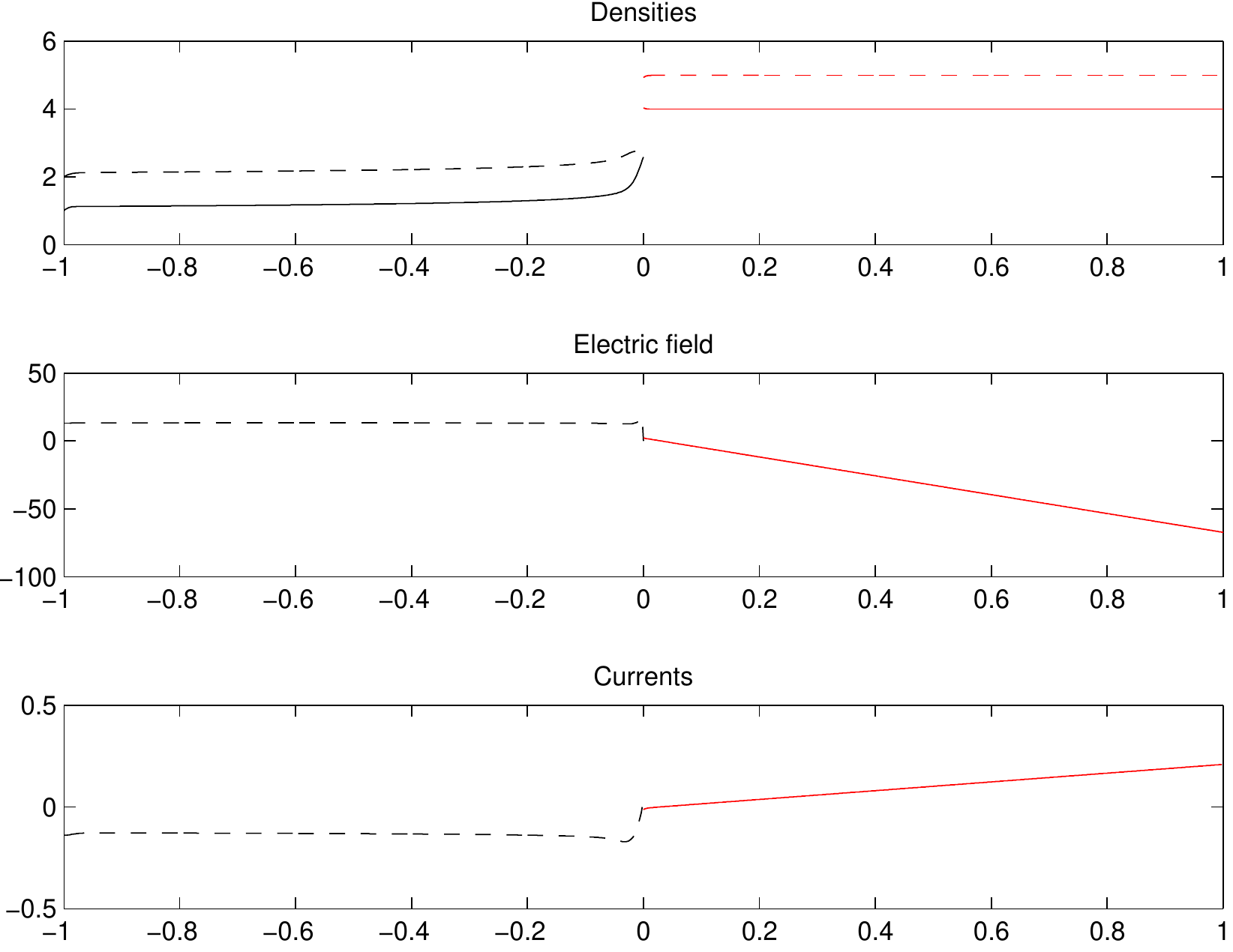} 
\includegraphics[height=0.3\textheight,width=0.4\textwidth]{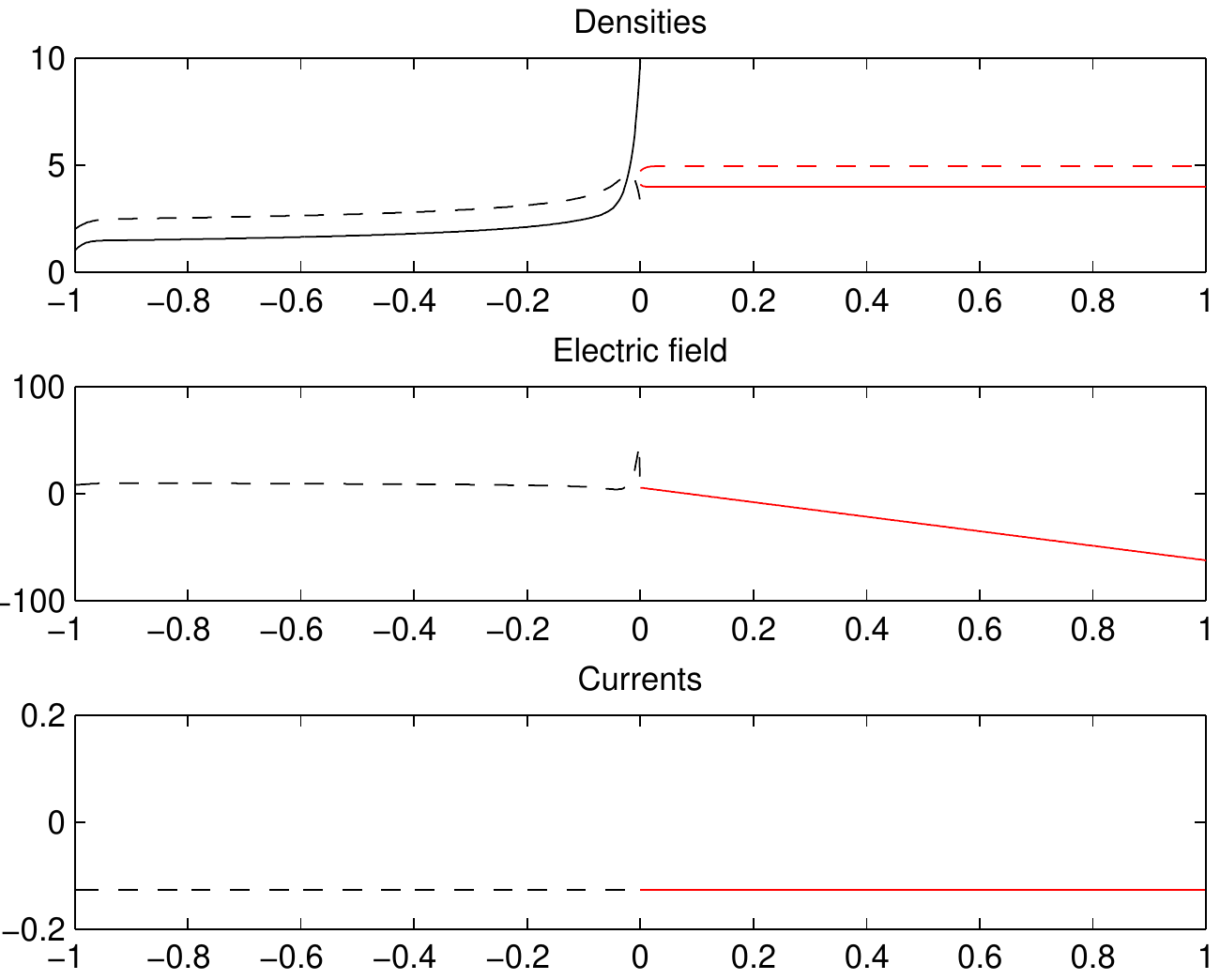}
\caption{Case I (c) in illuminated environment ($\gamma=1$). Shown are: charge densities (black solid for $\rho_p'$, black dashed for $\rho_n'$, red solid for $\rho_o'$ and red dashed for $\rho_r'$ respectively), electric field and flux distributions at time $t'=0.05$ (left column); and charge densities, electric field, and flux distributions at the stationary state (right column).}
\label{FIG:Ill I(c)}
\end{figure}

\paragraph{Case I (c).} There are a large number of physical parameters in the semiconductor-electrolyte system that control the dynamics of the system. To be specific, we show in this numerical simulation the impact of the relative dielectric constant in the electrolyte $\eps_r^\sfE$ on the system performance. The setup is exactly as in Case I (b) except that $\eps_r^\sfE=100$ now. We perform simulations in both the dark and illuminated environments. We plot in Fig.~\ref{FIG:Dark I(c)} the densities, the electric field, and the flux distributions at time $t'=0.05$ and the stationary state. The corresponding results for illuminated case are shown in Fig.~\ref{FIG:Ill I(c)}. It is not hard to observe the dramatic change in all the quantities shown after comparing Fig.~\ref{FIG:Dark I(c)} with Fig.~\ref{FIG:Dark I(b)}, and Fig.~\ref{FIG:Ill I(c)} with Fig.~\ref{FIG:Ill I(b)}. More detailed study of the parameter sensitivity problem will be reported elsewhere.

\subsection{Comparison with Schottky approximation}

In simulations done for practical applications, it is often assumed that the densities for the reductants and oxidants in the electrolyte are extremely high, so that reductant-oxidant dynamics changes very little compared to the electron-hole dynamics in the semiconductor; see, for instance, ~\cite{LaBa-JES76A,LaBa-JES76B}. In this case, it is usually common to completely fix the electrolyte system and only evolve the electron-hole system. This is done by treating the electrolyte system as a Schottky contact and thus replacing the semiconductor-electrolyte interface conditions with Robin-type of boundary conditions such as~\eqref{EQ:DDP BC Robin}; see more discussions in~\cite{LaBa-JES76A,LaBa-JES76B,Fawcett-Book04,JiWaSiKaLeXi-EES14,Memming-Book01}

Our previous simulations, such as those shown in Fig.~\ref{FIG:Dark I(a)} and Fig.~\ref{FIG:Ill I(a)}, indicate that such a simplification indeed can be accurate as the densities of the redox pair are roughly constant inside the electrolyte. We now present two simulations where we compare the simulations based on the our model of the whole system with the simulations based on Schottky approximation. Our simulations focus on Device II, which is significantly smaller than Device I. Simulations for Device I show similar behavior which we omit here to avoid repetition.

Due to the fact that the Schottky boundary condition contains no information on the parameters of the electrolyte system, it is impossible to make a direct quantitative comparison between the simulations. Instead, we will perform a comparison as follows. For each configuration of the electrolyte system, we select the parameters $\Phi_m$, $\chi$ and $E_g$ such that total potential on the Schottky contact, $\varphi_{Stky}'+\varphi_{app}'$, equals the applied potential on the counter electrode, $\varphi_{app}^{\sfA'}$.
\begin{figure}[ht]
\centering
\includegraphics[height=0.3\textheight,width=0.45\textwidth]{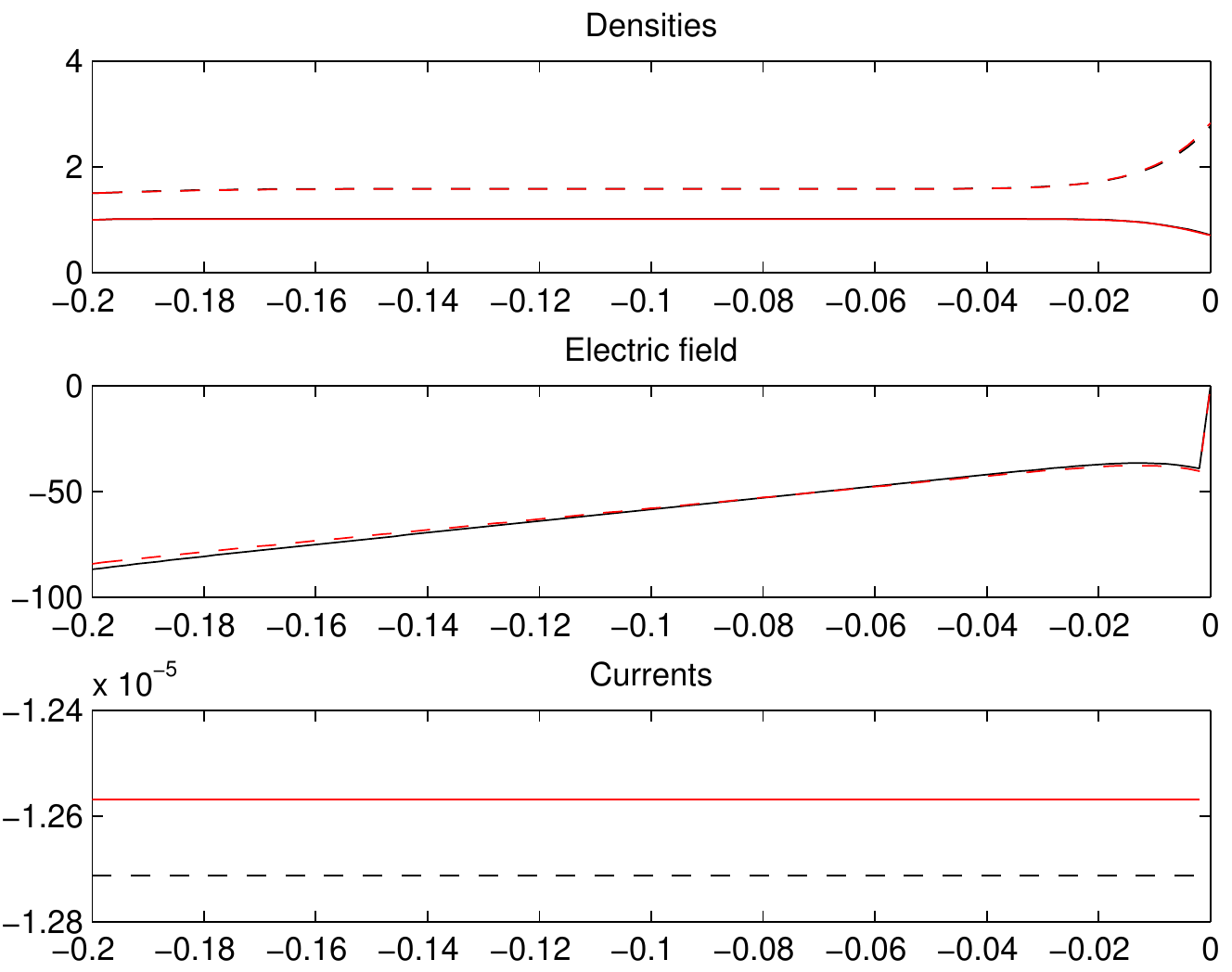}
\includegraphics[height=0.3\textheight,width=0.45\textwidth]{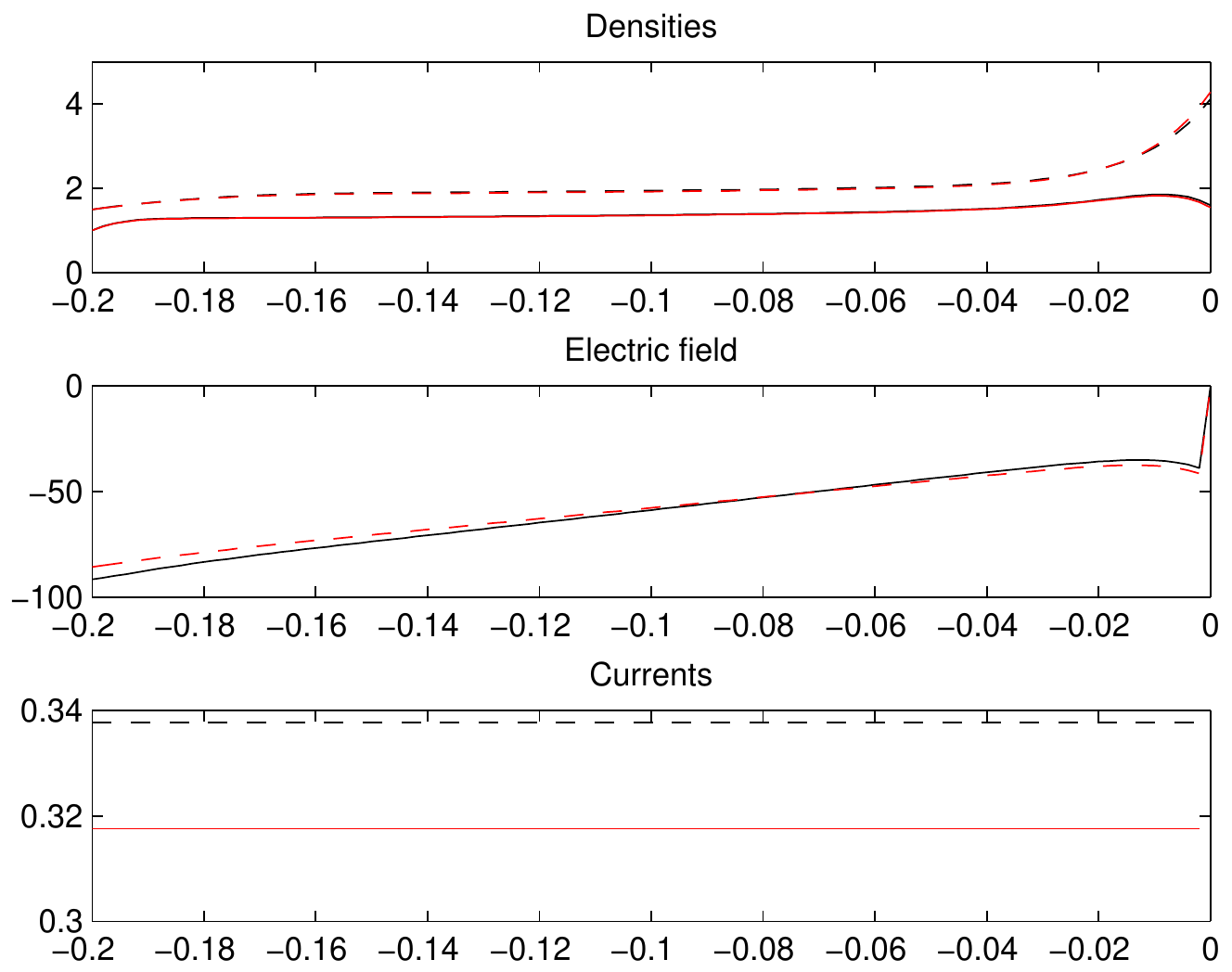}
\caption{Case II (a): comparison between simulations from whole-system (red) and Schottky approximation (black) for Device II in dark ($\gamma=0$, left column) and illuminated ($\gamma=1$, right column) environments.}
\label{FIG:II(a)}
\end{figure}

We performed two sets of simulations on Device II. In the simulations with Schottky approximation, the Schottky contact locates at $x'=0$. 

\paragraph{Case II (a).} In the first set of simulations, we compare the whole system simulation with the Schottky approximation when the bulk reductant-oxidant pair density are more than $20$ times higher than that of the electron-hole pair density. Precisely, the bulk densities for reductant and oxidant are respectively $\rho_r^{'\infty}=30.0$ and $\rho_o^{'\infty}=35.0$. The results in the semiconductor part are shown in Fig.~\ref{FIG:II(a)}. 
We observed that in this case, the simulations with Schottky approximation are sufficiently close to the simulations with the whole system. This is to say that, under such a case, replacing the electrolyte system with a Schottky contact provides a good approximation to the original system.

\begin{figure}[ht]
\centering
\includegraphics[height=0.3\textheight,width=0.45\textwidth]{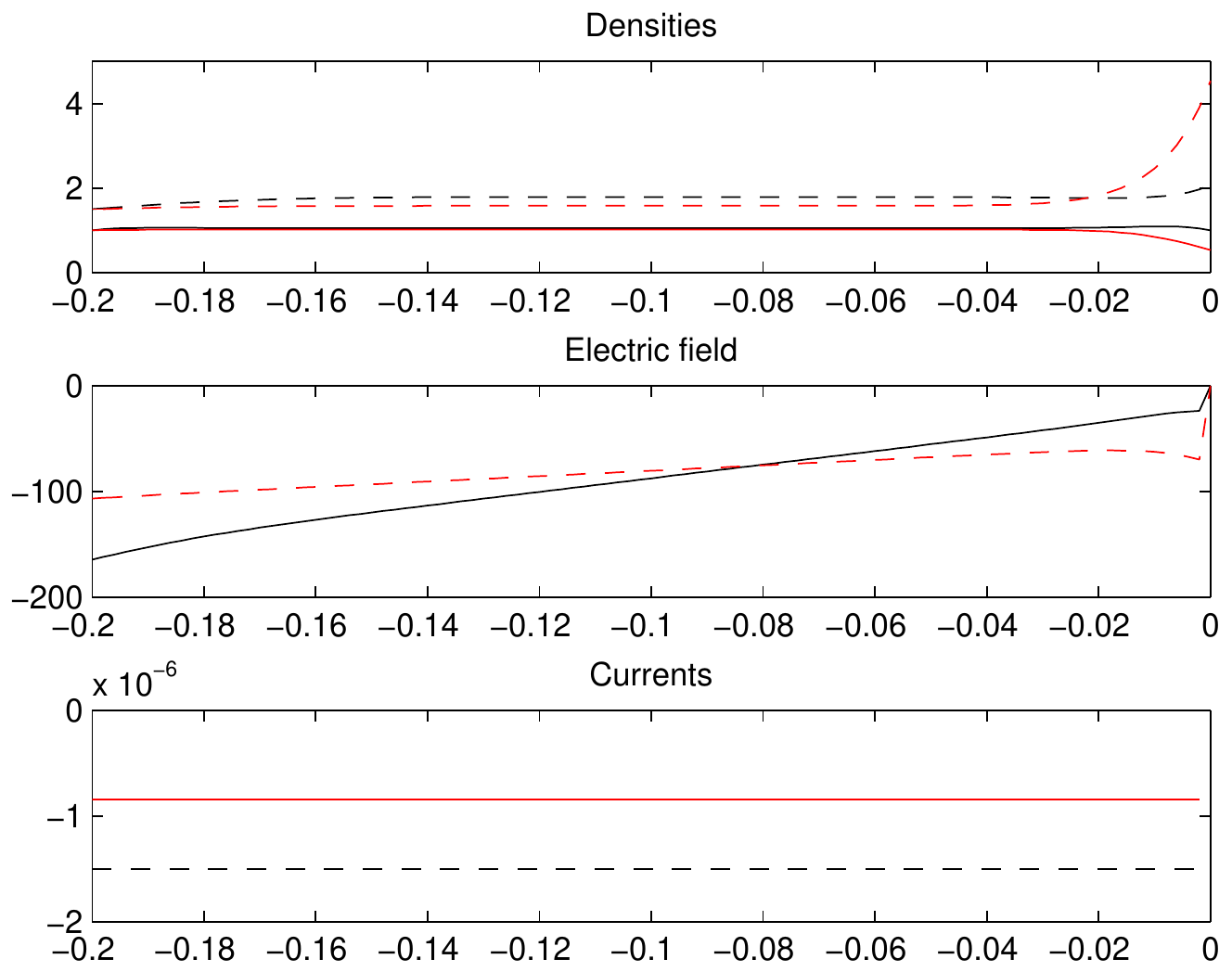}
\includegraphics[height=0.3\textheight,width=0.45\textwidth]{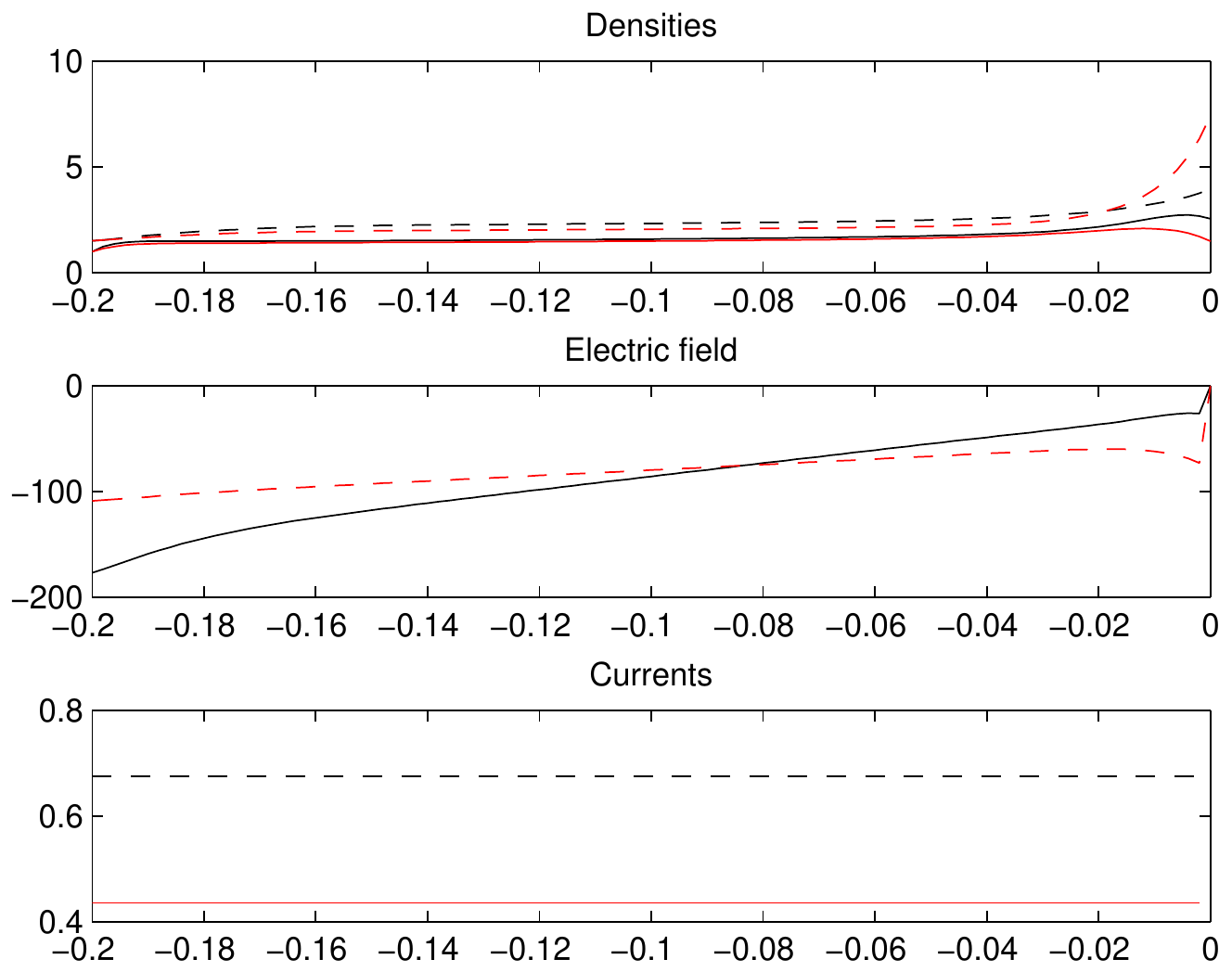}
\caption{Case II (b): comparison between simulations from whole-system (red) and Schottky approximation (black) for Device II in dark ($\gamma=0$, left column) and illuminated ($\gamma=1$, right column) environments.}
\label{FIG:II(b)}
\end{figure}

\paragraph{Case II (b).} In the second set of simulations, we compare the two simulations when the bulk reductant-oxidant pair density is much lower than that in Case II(a). Precisely, $\rho_r^{'\infty}  = 2.0$ and $\rho_o^{'\infty}  = 3.0$. The results for the semiconductor part are shown in Fig.~\ref{FIG:II(b)}.
While the qualitative behavior of the two systems is similar, the quantitative results are very different. We were not able to find a set of parameters in the boundary condition~\eqref{EQ:DDP BC Robin} that produced identical behavior of the semiconductor system. In this type of situations, simulation of the whole semiconductor-electrolyte system provides more accurate description of the physical process in the device.

The results above show that even though one can often replace the electrolyte system with a Schottky contact, modeling the whole semiconductor-electrolyte system offers more flexibilities in general. In the case when the bulk densities of the reductant and the oxidant are not extremely high compared to the densities of electrons and holes in the semiconductor, replacing the electrolyte system with a Schottky contact can cause large inaccuracy in predicted quantities.

\subsection{Voltage-current characteristics}

In the last set of numerical simulations, we study the voltage-current, or more precisely the voltage-flux (since all quantities are nondimensionalized), a characteristic of the semiconductor-electrolyte system. We again focus on simulations with Device II, although we observe similar results for Device I.

\begin{figure}[ht]
\centering
\includegraphics[height=0.2\textheight,width=0.45\textwidth]{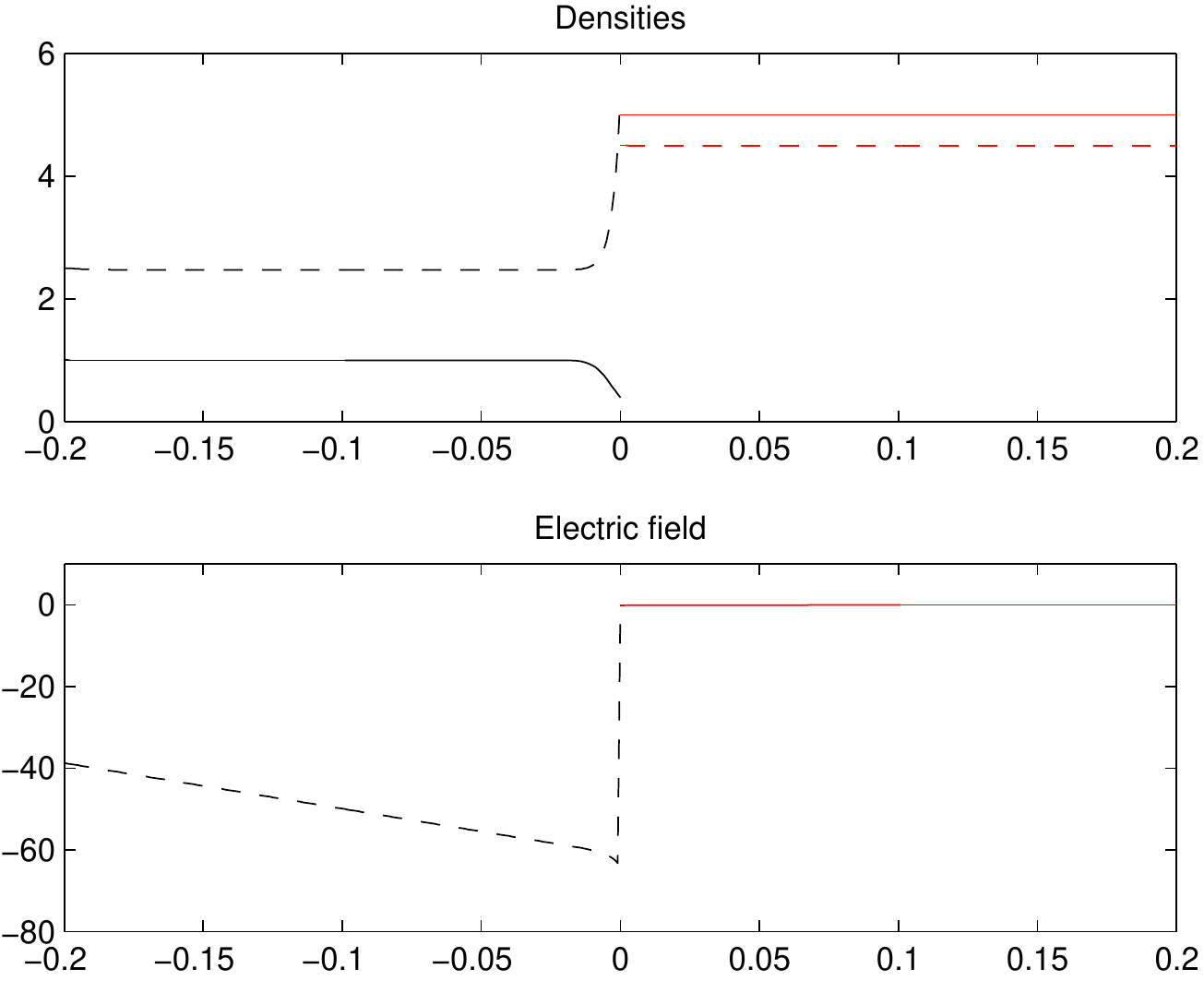}
\includegraphics[height=0.2\textheight,width=0.45\textwidth]{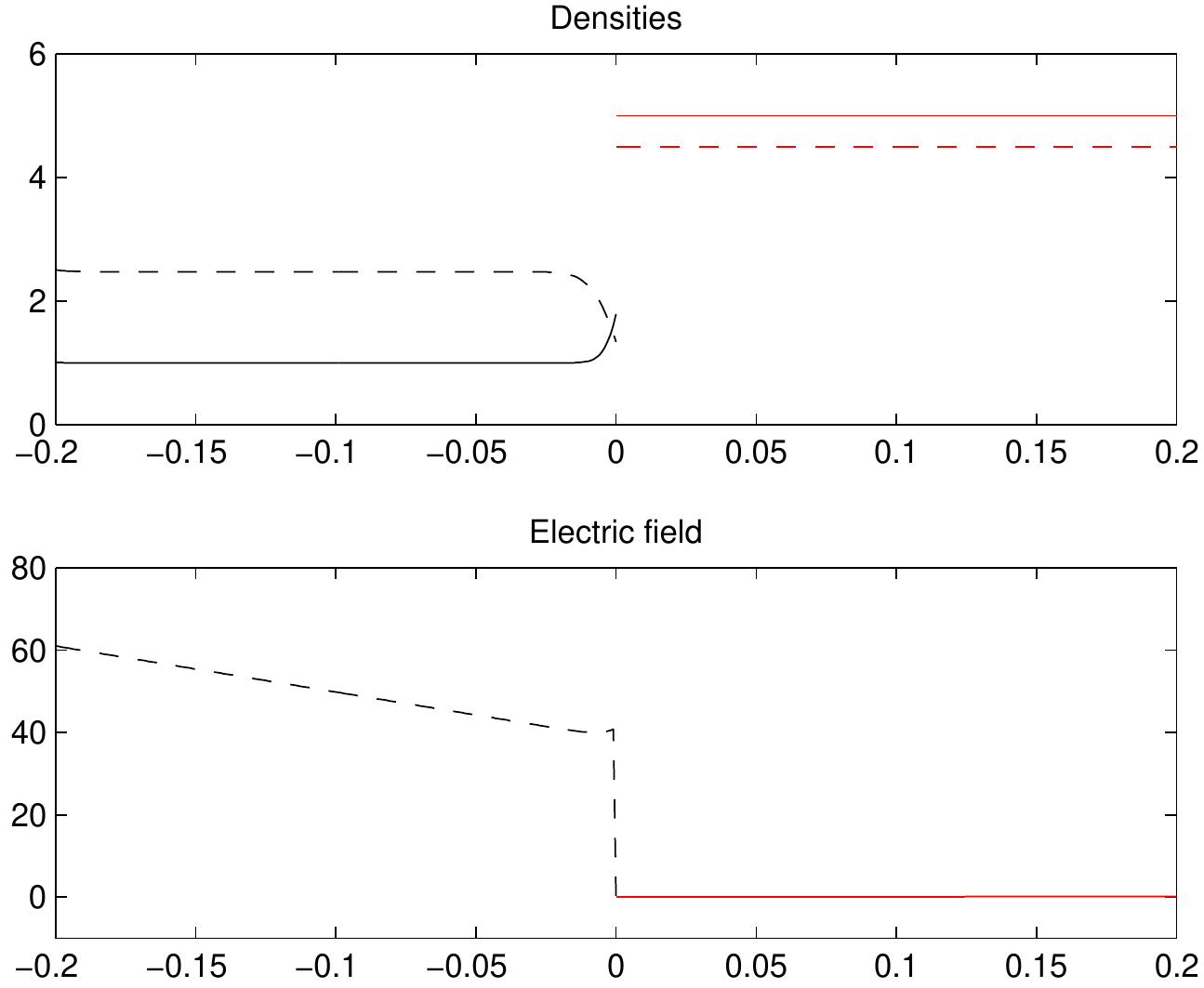}
\caption{Case II$^\prime$ (a)  in dark ($\gamma=0$) environment: comparison of charge densities (top) and the corresponding electric fields (bottom) in applied forward (left column) and reversed (right column) potential bias.}
\label{FIG:Dark IIp(a)}
\end{figure}

\begin{figure}[ht]
\centering
\includegraphics[height=0.2\textheight,width=0.6\textwidth]{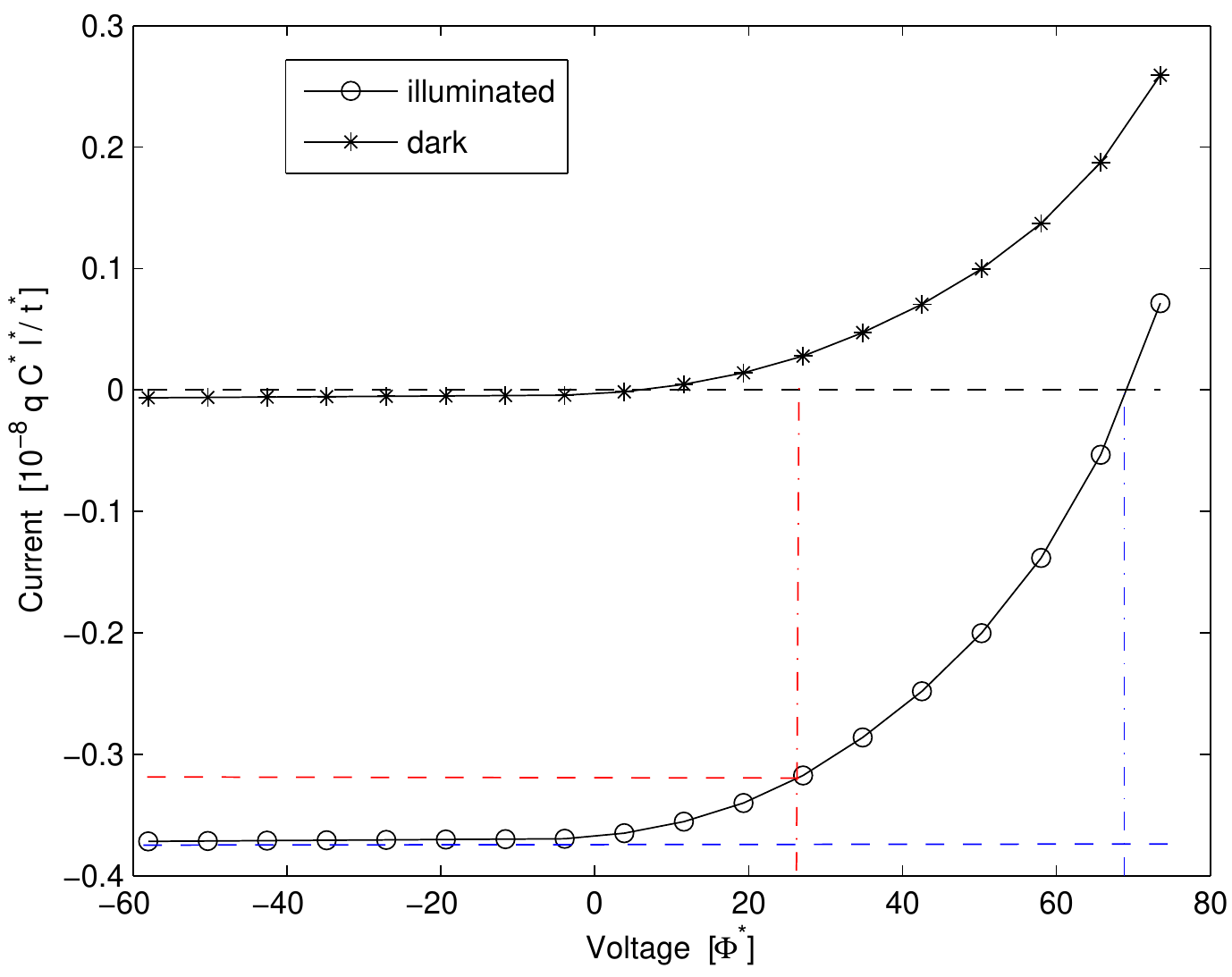}
\caption{Case II$^\prime$ (b): voltage-flux (I-V) curves for an $n$-type semiconductor-electrolyte system in dark (dotted line with circles) and illuminated (solid line with dots) environments. The intersections of the dash-dotted vertical lines with the x-axis show the maximum power voltage ($\Phi_{\rm mp}$, red) and open circuit voltage ($\Phi_{\rm oc}$, blue). The intersections of the dashed horizontal lines with the y-axis show the maximum power flux ($J_{\rm oc}$, red) and short circuit flux ($J_{\rm sc}$, blue).}
\label{FIG:IIp(b)}
\end{figure}

\paragraph{Case II$^\prime$ (a).} We first compare simulations with forward potential bias with simulations with reverse potential bias in a dark environment. We take a device with low bulk reductant-oxidant densities ($\rho_r^{'\infty}=4.5$, $\rho_o^{'\infty}=5.0$) and high relative dielectric constant in electrolyte ($\eps_r^\sfE=1000$). Shown on the left column of Fig.~\ref{FIG:Dark IIp(a)} are electron-hole and redox pair densities and the corresponding electric fields for applied forward potential bias of $35.0$. The same results for the applied reverse bias of $-35.0$ are presented in the right column of the same figure. Aside from the obvious differences in the densities and electric field distributions, the fluxes through the system with forward and reverse biases are very different, as can be seen later on Fig.~\ref{FIG:IIp(b)}.

\begin{figure}[ht]
\centering
\includegraphics[height=0.2\textheight,width=0.45\textwidth]{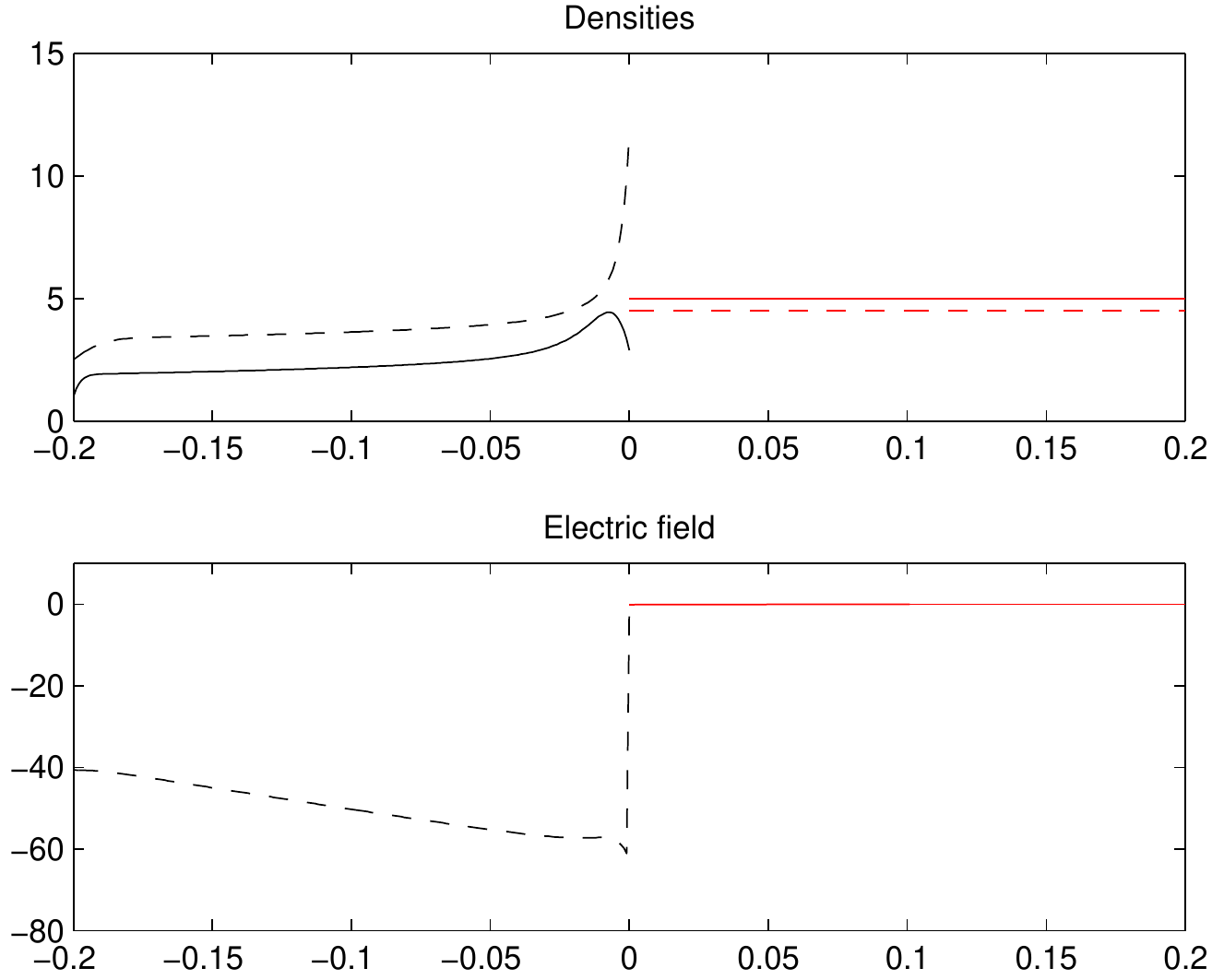}
\includegraphics[height=0.2\textheight,width=0.45\textwidth]{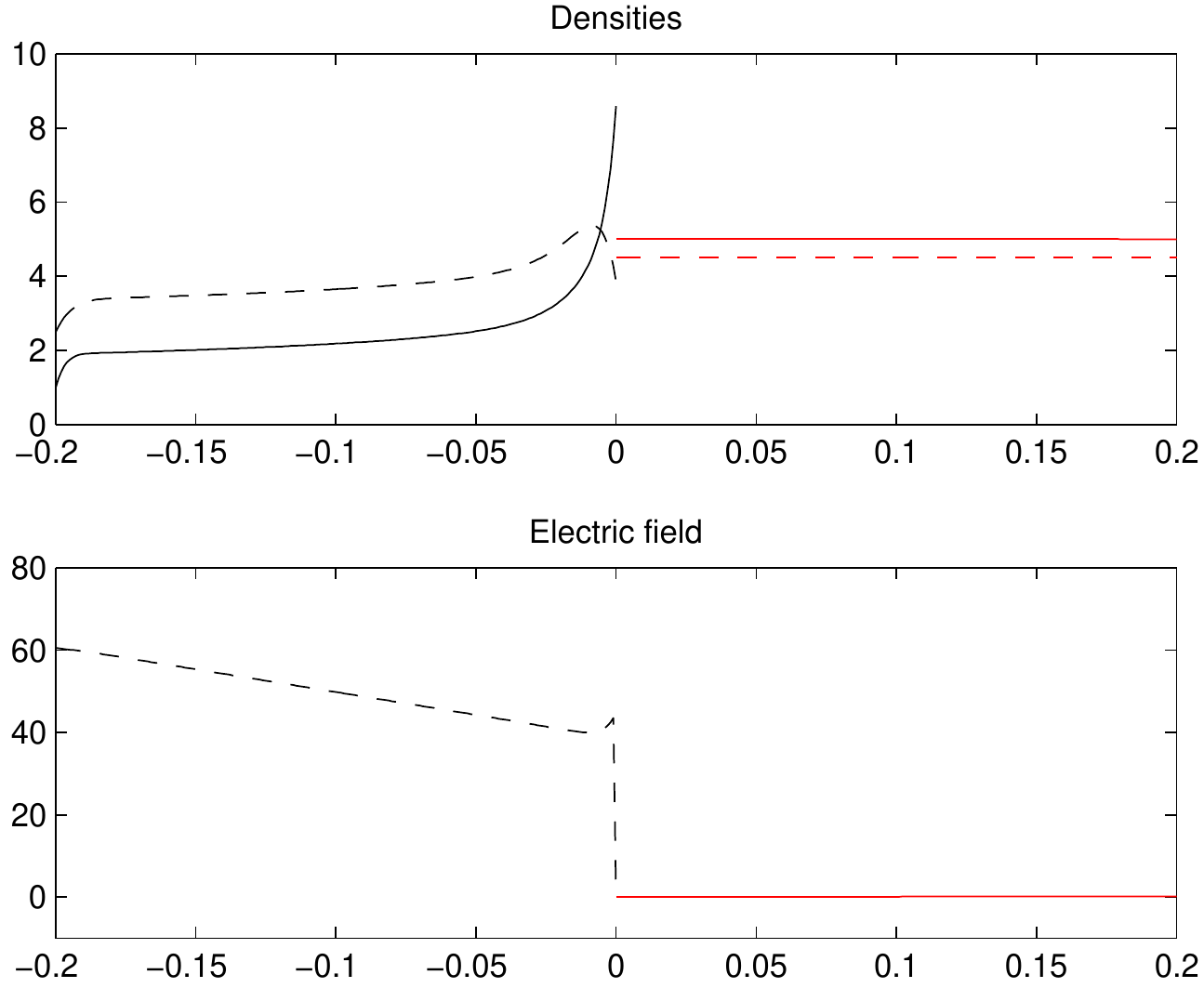}
\caption{Case II$^\prime$ (a) in illuminated ($\gamma=1$) environment: comparison of charge densities (top) and the corresponding electric fields (bottom) in applied forward (left) and reverse (right) potential bias.}
\label{FIG:Ill IIp(a)}
\end{figure}

The simulations are repeated in Fig.~\ref{FIG:Ill IIp(a)} in the illuminated environment. It is clear from the plots in Fig.~\ref{FIG:Dark IIp(a)} and Fig.~\ref{FIG:Ill IIp(a)} that illumination dramatically changes the distribution of charges (and thus the electric field) inside the device as we have seen in the previous simulations. 

\paragraph{Case II$^\prime$ (b).} We now attempt to explore the whole voltage-flux (I-V) characteristics of Device II. To do that, we run the simulations for various different applied potentials and compute the flux through the system under each applied potential. We plot the flux data as a function of the applied potential to obtain the I-V curve of the system. The parameters are taken as $\rho_r^{'\infty}=35$, $\rho_o^{'\infty}=30$), and $\eps_r^\sfE=1000$ to mimic those in realistic devices. 
Fig.~\ref{FIG:IIp(b)} shows the I-V curve obtained in both dark (line with stars) and illuminated (solid line with circles) environments. The applied potential bias lives in the range $[-58.0,\ 73.5]$. The intersections of the dash-dotted vertical lines with the $x$-axis (horizontal) show the maximum power voltage ($\Phi_{\rm mp}$, red) and open circuit voltage ($\Phi_{\rm oc}$, blue). The intersections of the dashed horizontal lines with the y-axis show the maximum power flux ($J_{\rm oc}$, red) and short circuit flux ($J_{\rm sc}$, blue).

\section{Conclusion and remarks}
\label{SEC:Concl}

We have considered in this paper the mathematical modeling of semiconductor-electrolyte systems for applications in liquid-junction solar cells. We presented a complete mathematical model, a set of nonlinear partial differential equations with reactive interface conditions, for the simulation of such systems. Our model consists of a reaction-drift-diffusion-Poisson system that models the transport of electron-hole pairs in the semiconductor and an equivalent system that describes the transport of reductant-oxidant pairs in the electrolyte. The coupling of the two systems on the semiconductor-electrolyte interface is modeled with a set of reaction and flux transfer types of interface conditions. We presented numerical procedures to solve both the time-dependent and stationary problems, for instance with Gummel-Schwarz double iteration. Some numerical simulations for one-dimensional devices were presented to illustrate the behavior of these devices.

Past study on the semiconductor-electrolyte system usually completely neglected the charge transfer processes in the electrolyte~\cite{LaBa-JES76A,LaBa-JES76B,Fawcett-Book04,JiWaSiKaLeXi-EES14,Memming-Book01}. The mathematical models developed thus only cover the semiconductor part with the interface effect modeled by a Robin type of boundary condition. The rationale behind this simplification is the belief that the density of the reductant-oxidant pairs is so high compared to the density of the electron-hole pairs in the semiconductor such that the density of the redox pair would not be perturbed by the charge transfer process through the interface. While this might be a valid approximation in certain cases, it can certainly go wrong in other cases. For instance, it is generally observed that due to the strong electric field at the interface, there are dramatic change in the density of charges (both the electron-hole and the redox pairs) near the interface. What we have presented is, to our knowledge, the first complete mathematical model for semiconductor-electrolyte solar cell systems that would allow us to accurately study the charge transfer process through the interface. 

There are a variety of problems related to the model that we have constructed in this work that deserve thorough investigations. On the mathematical side, a detailed mathematical analysis on the well-posedness of the system is necessary. On the computational side, more detailed numerical analysis of the model, including convergence of the Gummel-Schwarz iteration, efficient high-order discretization, and fast solution techniques, has to be studied. On the application side, it is important to calibrate the model parameters with experimental data collected from real semiconductor-electrolyte solar cells. Once the aforementioned issues are addressed, we can use the model to help design more efficient solar cells by, for instance, optimizing the various model parameters. We are currently investigating several of these issues.

\section*{Acknowledgements}

We would like to thank the anonymous referees for their constructive comments that helped us improve the quality of this work. We also benefited greatly from fruitful discussions with Professor Allen J. Bard and Charles B. Mullins (both affiliated with the Department of Chemistry at the authors' institution). YH, IMG and HCL are partially supported by the National Science Foundation through grants CHE 0934450 and DMS-0807712. KR is partially supported by NSF grants DMS-0914825 and DMS-1321018.



\end{document}